\documentclass[11pt, a4paper, english]{article}

\usepackage{placeins}                 
\usepackage[english]{babel}
\usepackage{graphicx, wrapfig, subfig} 
\usepackage{amsmath, amsfonts} 
\usepackage{amsthm} 
\usepackage{amssymb} 
\usepackage{mathrsfs} 
\usepackage[retainorgcmds]{IEEEtrantools} 
\usepackage{xcolor} 
\usepackage{adjustbox}
\usepackage{algpseudocode} 
\usepackage{algorithm} 
\graphicspath{{./Figures/}} 
\theoremstyle{remark}
\newtheorem{remark}{Remark}
\newtheorem{simplif}{Simplification}

\newcommand{\norm}[1]{\left\lVert#1\right\rVert}

\newcommand{\BR}{\mathbb{R}}

\renewcommand{\d}{\;\mathrm{d}}
\newcommand{\N}{\mathcal{N}}
\newcommand{\M}{\mathcal{M}}
\newcommand{\J}{\mathcal{J}}
\renewcommand{\P}{\mathcal{P}}
\newcommand{\V}{\mathcal{V}}

\newcommand{\D}{\mathcal{D}}
\renewcommand{\H}{\mathrm{H}}
\newcommand{\K}{\mathcal{K}}

\title{Isogeometric Methods for Free Boundary Problems}
\author{M. Montardini\textsuperscript{1}, F. Remonato\textsuperscript{1,2}, G. Sangalli\textsuperscript{1,3} \\[10pt]
		\small \textsuperscript{1} Department of Mathematics, University of Pavia, Pavia, Italy\\
		\small \textsuperscript{2} Department of Mathematical Sciences, NTNU, Trondheim, Norway\\[5pt]
		\small \textsuperscript{3} IMATI-CNR ``E. Magenes'', Pavia, Italy\\
		\small monica.montardini01@universitadipavia.it, \\ \small filippo.remonato@ntnu.no, \\ \small giancarlo.sangalli@unipv.it }

\date{}

\begin{document}

\maketitle

\begin{abstract}
	We present in detail three different quasi-Newton isogeometric algorithms for the treatment of free boundary problems.
	Two algorithms are based on standard Galerkin formulations, while the third is a fully-collocated scheme.
	With respect to standard approaches, isogeometric analysis enables the accurate description of curved geometries, and is thus particularly suitable for free boundary numerical simulation.
	We apply the algorithms and compare their performances to several benchmark tests, considering both Dirichlet and periodic boundary conditions.
	In this context, iogeometric collocation turns out to be robust and computationally more efficient than Galerkin.
	Our results constitute a starting point of an in-depth analysis of the Euler equations for incompressible fluids.
\end{abstract}


\section{Introduction}

This work focuses on the isogeometric analysis (IGA) of free boundary problems.
IGA, first presented in \cite{hughes2005isogeometric}, is a recent extension of the standard finite element method where the unknown solution of the partial differential equation is approximated by the same functions that are adopted in computer-aided design for the parametrization of the problem domain.
These functions are typically splines and extensions, such as non-uniform rational B-splines (NURBS).
We refer to the monograph \cite{Cottrell2009iat} for a detailed description of this approach.

In this work we present three general free boundary algorithms.
The first algorithm is an extension to IGA of the finite elements approach of \cite{Karkkainen1999fss, Karkkainen2004sca}.
Since the finite element basis produces meshes with straight edges, the authors needed a workaround to approximate the curvature of the boundary; in the new IGA framework this can be avoided thanks to the natural description of curved geometries through spline functions.
IGA of free boundary problems was already considered in \cite{VanderZee2010goe, Vanderzee2013snm}; our second algorithm uses and extends these approaches to problems with periodic conditions. 
Our third and most efficient scheme uses instead an isogeometric variational collocation approach based on the superconvergent points presented in \cite{gomez2016variational,Montardini2016ooi}.
The choice of applying an IGA collocation method is a novelty in this setting and, moreover, allows for a fast computation of the solution.
While speed is marginally important in the benchmarks considered in this work, it becomes a major concern when one needs to address more complicated problems.

All the algorithms are based on shape calculus techniques, see for example \cite{Delfour2011sag, Sokolowski1992its}.
This results in the three algorithms being of quasi-Newton type, achieving superlinear convergence.

Our interest in free boundary problems is motivated by a separate analysis, in progress at the time of writing, of the periodic solutions of the Euler equations describing the flow of an incompressible fluid over a rigid bottom.
The analytical literature on this problem is quite extensive, with results regarding irrotational flows \cite{Groves2004sww}, the limiting Stokes waves \cite{Toland1996sw}, or waves on a rotational current containing one or multiple critical layers \cite{Ehrnstrom2011sww, Wahlen2009sww}.
The numerical experiments so far have used finite differences methods \cite{Darlymple1977anm}, boundary-integral formulations \cite{Simmen1985sdw}, or finite elements \cite{Rycroft2013cot}.
Several other examples and numerical experiments, also based on boundary formulations, can additionally be found in \cite{Vadenbroeck2010gcf}.

This paper is organised as follows: In Section \ref{sec:FBP} we describe the details of free boundary problem, and present two weak formulations that will constitute our starting point for the algorithms.
In Section \ref{sec:linearisation} we first introduce the necessary shape calculus tools, and then proceed to linearise the aforementioned weak forms.
This will produce the correct formulations on which to base our quasi-Newton steps.
Section \ref{sec:schemes} describes the discrete spaces used in the numerical schemes along with the structure of the algorithms.
Finally, Section \ref{sec:results} presents the numerical benchmarks and the results we obtained.
We summarise the results and draw our conclusions in Section \ref{sec:conclusions}.


\section{Free Boundary Problem}
\label{sec:FBP}

\begin{figure}
	\centering
	
	\includegraphics[width=0.7\linewidth]{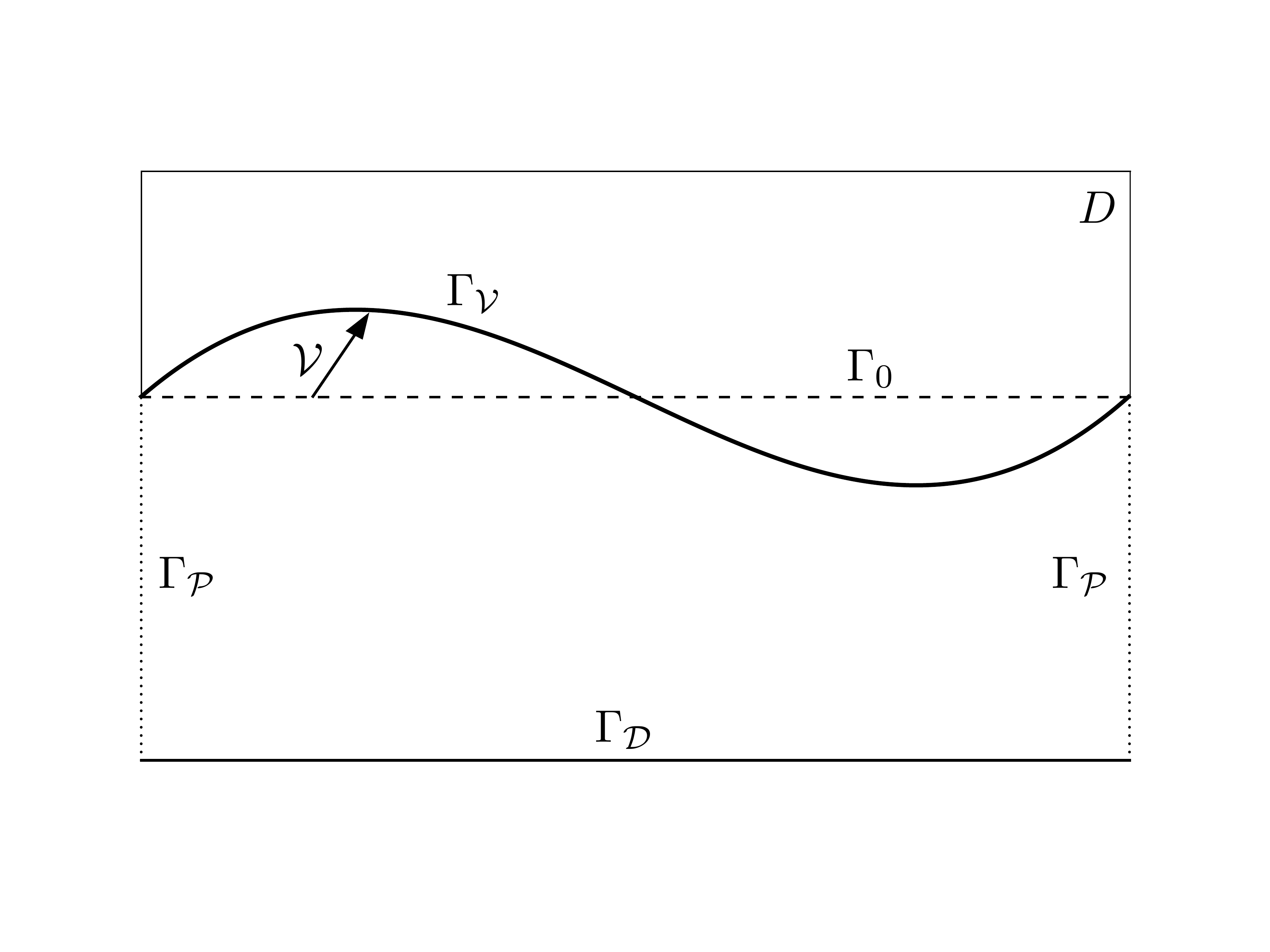}
	
	\caption{\label{fig:Setting}The setting of our problem. The vector field $\V$ deforms the reference free boundary $\Gamma_0$ (dashed line) into the free boundary $\Gamma_\V$ (thick solid line). The vertical dotted lines represent the periodic boundary $\Gamma_\P$, while the thin solid line represents the fixed flat bottom boundary $\Gamma_\D$. The physical domain and its deformations are contained in a larger rectangle $D$.}
\end{figure}

Let $\Omega_0$ be a domain used as reference configuration with $\partial\Omega_0 = \Gamma_\D \cup \Gamma_\P \cup \Gamma_0$;  $\Gamma_\D$ being the (fixed) bottom boundary with Dirichlet data, $\Gamma_\P$ the (fixed) vertical boundary with periodic conditions, and $\Gamma_0$ the (free) upper part of the boundary.
Moreover, let $D$ be a rectangle with basis $\Gamma_D$, containing $\Omega_0$ and all its possible deformations. 
For $M$ a domain and $\Gamma$ a curve, we denote with $C^{k,\lambda}(M,\mathbb{R}^2)$ the space of $(k, \lambda)-$H\"older continuous functions defined on $M$ with values in $\mathbb{R}^2$ and by $C^{k,\lambda}_{0}(\Gamma,\mathbb{R}^2)$ the subspace of $C^{k,\lambda}(\Lambda;\mathbb{R}^2)$  with compact support, in particular vanishing at the two extremes of the curve. 
Then, the set of admissible vector fields acting on the reference domain is defined as $\Theta=\{ \V \in C^{0,1}(D,\mathbb{R}^2) \cap C^{1,1}_{0}(\Gamma_0,\mathbb{R}^2) \,|\, \V=0 \ \text{on} \ \Gamma_\D \text{ and } \V(\cdot,y) \text{ periodic}\}$. 
We encode the deformation of the upper part of the boundary, $\Gamma_0$, as the action of a vector field $\V \in \Theta$ such that the deformed domain is smooth enough, does not have self intersections and does not touch the bottom $\Gamma_D$.
For this reason we denote the deformed free boundary with $\Gamma_\V = \{ x\in \BR^2 \,|\, x = x_0 + \V(x_0), \, x_0 \in \Gamma_0\}$.
Analogously, $\Omega_\V$ will denote the physical domain with boundary $\partial\Omega_\V = \Gamma_\D \cup \Gamma_\P \cup \Gamma_\V$; see Figure \ref{fig:Setting} for a representation of this setting.
We remark that $\Gamma_0$ is in general not flat.

The Bernoulli-type free boundary problem (FBP) we are interested in can then be posed as searching for a pair $(u,\V)$, both periodic in the $x$-direction, such that

\begin{IEEEeqnarray}{rCls}
	\IEEEyesnumber \label{eq:FBP} \IEEEyessubnumber*
	-\Delta u &=& f &\quad in $\Omega_\V$ \label{eq:FBP1}\\
	u &=& h &\quad on $\Gamma_\V \cup \Gamma_\D$ \label{eq:FBP2}\\
	\partial_n u &=& g &\quad on $\Gamma_\V$ \label{eq:FBP3}
\end{IEEEeqnarray} 
where $\partial_n u = \nabla u \cdot n$ is the outward normal derivative of $u$.
The functions $f$, $h$, and $g$ are defined in $D$ and are compatible with the periodicity requirement.
We will consider $h$ and $g$ continuous, with $g$ strictly positive and bounded away from zero\footnote{The strict positivity is not strictly necessary: If $g<0$ one could, for instance, keep track of the sign of $g$ in the numerical method directly. However, $g$ has to have a definite sign everywhere on $\Gamma_\V$.}.
\begin{remark}
	The analytical treatment of the problem with periodic boundary conditions does not differ much from the case with pure Dirichlet conditions, which we also consider in our numerical benchmarks. 	
\end{remark}

\subsection{Weak Formulation}
To obtain a formulation of \eqref{eq:FBP} suitable for a numerical scheme we first follow the steps presented in \cite{Karkkainen1999fss}.
This approach leads to two distinct, coupled weak forms.
Given the space $H^1_{per}(\Omega_\V) = \{ u \in H^1(\Omega_\V) \,|\, u(\cdot , y) \text{ periodic}\}$, for a known function $r$ periodic in the $x$-direction we define the space
\begin{IEEEeqnarray*}{rCl}
	H^1_{r, \Gamma_\D}(\Omega_\V)& = & \{ \varphi \in H^1_{per}(\Omega_\V) \,|\, \varphi = r \text{ on } \Gamma_\D\}.
\end{IEEEeqnarray*}

The first weak form is then obtained using \eqref{eq:FBP1}, \eqref{eq:FBP3}, and the part of \eqref{eq:FBP2} pertaining to $\Gamma_\D$. 
We select test functions $\varphi \in H^1_{0, \Gamma_\D}(\Omega_\V)$ and apply Green's formula once to obtain
\begin{equation}
	\label{eq:startingWeak1}
	\int_{\Omega_\V} \nabla u \cdot \nabla \varphi \d\Omega - \int_{\Gamma_\V} g\,\varphi \d\Gamma = \int_{\Omega_V} f\, \varphi \d\Omega.
\end{equation}
Using the part of \eqref{eq:FBP2} on $\Gamma_\V$ we employ test functions $v \in H^1_{per}(\Gamma_\V)$ and write the second weak form simply as
\begin{equation}
	\label{eq:startingWeak2}
	\int_{\Gamma_\V} \!uv \d\Gamma = \int_{\Gamma_\V} \!hv \d\Gamma.
\end{equation}
We select the trial function space by requiring $ u \in H^1_{h, \Gamma_\D}(\Omega_\V)$, thereby strongly imposing the Dirichlet boundary conditions on $\Gamma_\D$. 
This leads to the definition of two linear forms:
\begin{IEEEeqnarray}{rCl}
	\label{eq:M1}
	\M_1(u,\V; \varphi) &=& \int_{\Omega_\V} \nabla u \cdot \nabla \varphi \d\Omega - \int_{\Gamma_\V} g\,\varphi \d\Gamma - \int_{\Omega_\V} f\, \varphi \d\Omega, \\[5pt]
	\label{eq:M2}
	\M_2(u,\V; v) &=& \int_{\Gamma_\V} \!uv \d\Gamma - \int_{\Gamma_\V} \!hv \d\Gamma.
\end{IEEEeqnarray}
Thus, with this approach the problem is defined as: Search  for $(u,\V) \in H^1_{h, \Gamma_\D}(\Omega_\V) \times \Theta$ such that
\begin{IEEEeqnarray*}{rCl}
	\M_1(u,\V; \varphi) & = & 0,\\
	\M_2(u,\V; v) & = & 0,
\end{IEEEeqnarray*}
for all test functions $(\varphi, v )\in  H^1_{0, \Gamma_\D}(\Omega_\V)\times H^1_{per}(\Gamma_\V)$.

\subsection{Very-Weak Formulation}
We now  follow the approach of \cite{Vanderzee2013snm}.
The main difference from the previous formulation is that we write a single very-weak formulation containing information from all boundary conditions.

Considering the subspace $H^2_{0, \Gamma_\D}(\Omega_\V) = \{ \varphi \in H^1_{0, \Gamma_\D}(\Omega_\V) \,|\, \varphi \in H^2(\Omega_\V)\}$, we multiply \eqref{eq:FBP1} by a test function $\varphi \in H^2_{0, \Gamma_\D}(\Omega_\V)$; integrating by parts twice leads to
\begin{equation}
	\label{eq:startingVeryWeak}
	-\int_{\Omega_\V} (u - h)\, \Delta \varphi \d\Omega  + \int_{\Omega_\V} \nabla h \cdot \nabla \varphi \d\Omega = \int_{\Omega_\V} f\,\varphi \d\Omega + \int_{\Gamma_\V} \varphi\,g \d\Gamma,
\end{equation}
which we demand to be satisfied for all $\varphi \in H^2_{0, \Gamma_\D}(\Omega_\V)$.
In view of the above formulation we can then select the trial function space simply as $H^1_{per}(\Omega_\V)$.
The Dirichlet boundary conditions are therefore all imposed weakly.

From Equation \eqref{eq:startingVeryWeak} we define the linear form
\begin{IEEEeqnarray}{rCl}
	\label{eq:N}
	\N(u,\V; \varphi) &=& -\int_{\Omega_\V} (u - h)\, \Delta \varphi \d\Omega \, + \int_{\Omega_\V} \nabla h  \cdot \nabla \varphi \d\Omega \IEEEnonumber\\[5pt]
	&& -\int_{\Omega_\V} f\,\varphi \d\Omega \, - \int_{\Gamma_\V} \varphi\,g \d\Gamma.
\end{IEEEeqnarray}

Thus, with this approach the problem is defined as: Search  for $(u,\V) \in H^1_{per}(\Omega_\V) \times \Theta$ such that
\begin{equation*}
\N(u, \V; \varphi) = 0
\end{equation*}
for all test functions $\varphi \in H^2_{0, \Gamma_\D}(\Omega_\V)$.

Note that this very-weak formulation cannot be used directly to implement a numerical scheme, as the trial and test spaces are unbalanced.


\section{Linearising the FBP}
\label{sec:linearisation}

We now proceed in deriving a quasi-Newton algorithm to solve the free boundary problem.
The dependence on the domain's geometry is handled through shape calculus techniques to express the derivatives with respect to the vector field $\V$.

\subsection{Shape Derivatives}

Here we briefly state the shape calculus results we will need for the linearisation.
An in-depth analysis of the assumptions and regularity requirements can be found in the original work by Delfour, Zol\'esio, and Sokolowski \cite{Delfour2011sag, Sokolowski1992its}.
An overview of shape calculus presented with a more modern approach can also be found in \cite{VanderZee2010goe}.

Let $\mathcal{O}$ be a family of admissible (smooth enough) domains; a functional $\J$ is called a \emph{shape functional} if $\J : \mathcal{O} \rightarrow \BR$.
Note therefore that for a fixed function $u$ and test functions $\varphi$ and $v$, the maps defined by the linear forms introduced earlier are shape functionals provided we identify each element $\V \in \Theta$ with the domain $\Omega_\V$ in which $\Omega_0$ is deformed by the action of $\V$.

In the particular case of a domain functional  $\J(\V)=\int_{\Omega_\V}\psi \d\Omega$ and a boundary functional $\mathcal{F}(\V)=\int_{\Gamma_\V} \phi \d\Gamma$, with $\psi$ and  $\phi$ smooth functions in $\mathbb{R}^2$ independent of $\V$, the shape derivatives of $\J$ and $\mathcal{F}$ are described by the following \emph{Hadamard formulas}:
\begin{IEEEeqnarray*}{rCl}
	\IEEEyesnumber\label{eq:hadamard} \IEEEyessubnumber*
	\langle \, \partial_\V \J(\V), \, \delta \V \,\rangle &=& \int_{\Gamma_\V} \psi \, \delta \V \cdot n \d\Gamma \label{eq:hadamard_domain}\\
	\langle \, \partial_\V \mathcal{F}(\V), \, \delta \V \,\rangle &=& \int_{\Gamma_\V} \left( \partial_n \phi + \H \phi \right)  \delta \V \cdot n \d\Gamma \label{eq:hadamard_bdy}
\end{IEEEeqnarray*}
where $\delta\V\in\Theta$ is a perturbation of the vector field, $\H$ is the \emph{signed (additive) curvature} of $\Gamma_\V$ and $n$ is the normal vector pointing outward. 
In particular, considering a parametrization of the free boundary $\Gamma_{\V}$ defined as $\gamma(t) = (t, y(t))$, then
\begin{equation*}
\H = -\frac{y''}{\left[1+(y')^2\right]^{3/2}}.
\end{equation*}

\subsection{Linearisation of the weak formulation}
Let us first consider the linear forms \eqref{eq:M1} and \eqref{eq:M2}.
We want to linearise  $\M_1$ and $\M_2$ with respect to $u$ and $\V$ at an arbitrary approximated solution $(u^*, \V^*)\in H^{1}_{h,\Gamma_\D}(\Omega_{\V^*}) \times \Theta$.

Since the dependence of $\M_1$ and $\M_2$ on $u$ is affine, their G\^{a}teaux derivatives with respect to $u$ in the direction $\delta u \in H^1_{0,\Gamma_\D}(\Omega_{\V^*})$  are simply given by:
\begin{IEEEeqnarray*}{rCl}
	\IEEEyesnumber\label{eq:weakShape_u} \IEEEyessubnumber*
	\langle \, \partial_u  \M_1[u^*, \V^*; \varphi ] , \delta u \,\rangle &=&\int_{\Omega_{\V^*}}\nabla\delta u\cdot \nabla\varphi \d\Omega \label{eq:weakShape_u1} \\
	\langle \, \partial_u  \M_2[u^*, \V^*; v ] , \delta u \,\rangle &=& \int_{\Gamma_{\V^*}}\delta u \, v \d \Gamma. \label{eq:weakShape_u2} 
\end{IEEEeqnarray*}
The linearisation with respect to the vector field $ \V$ in the direction $\delta\V\in\Theta$ is performed using the Hadamard formulas; we obtain: 
\begin{IEEEeqnarray*}{rCl}
	\IEEEyesnumber\label{eq:weakShape_theta} \IEEEyessubnumber*
	\langle \, \partial_\V \M_1[u^*, \V^*; \varphi ] , \delta \V \,\rangle &= & \int_{\Gamma_{\V^*}}\nabla u^* \!\cdot \nabla\varphi \; \delta\V\cdot n \d\Gamma \IEEEnonumber \\
	& & -\int_{\Gamma_{\V^*}} \left[\K_\H \varphi \, + g\,\partial_n \varphi \right]\, \delta \V \cdot n \d\Gamma \label{eq:weakShape_V1}\\[5pt]
	\langle \, \partial_\V \M_2[u^*, \V^*; v ] , \delta \V \,\rangle &=& \int_{\Gamma_{\V^*}} \mkern-5mu \left(\partial_nu^*- \partial_n h + \H (u^* - h) \right)  v\; \delta\V \cdot n \d\Gamma \IEEEeqnarraynumspace \IEEEnonumber\\
	& & + \int_{\Gamma_{\V^*}} \mkern-5mu (u^* - h) \, \partial_n v \, \delta\V \cdot n \d\Gamma \label{eq:weakShape_V2}
\end{IEEEeqnarray*}
where $\K_\H= \partial_n g + \H g +f$, and $\H$ is the curvature of $\Gamma_{\V^*}$. 

A Newton step at the point $(u^*, \V^*)$ has then the following structure: Search for $\delta u \in H^1_{0,\Gamma_\D}(\Omega_{\V^*})$ and $\delta \V \in \Theta$ such that  
\begin{IEEEeqnarray*}{rCCCl}
	\IEEEyesnumber\label{eq:Newton} \IEEEyessubnumber*
	\langle \, \partial_u  \M_1[u^*, \V^*; \varphi ] , \delta u \,\rangle & + & \langle \, \partial_\V \M_1[u^*, \V^*; \varphi ] , \delta \V \,\rangle  &=& -\M_1(u^*, \V^*; \varphi) \IEEEeqnarraynumspace \label{eq:Newton1}\\
	\langle \, \partial_u  \M_2[u^*, \V^*; v ] , \delta u \,\rangle & + & \langle \, \partial_\V \M_2[u^*, \V^*; v ] , \delta \V \,\rangle  &=& -\M_2(u^*, \V^*; v) \label{eq:Newton2}
\end{IEEEeqnarray*}
for all $(\varphi, v )\in  H^1_{0, \Gamma_\D}(\Omega_{\V^*})\times H^1_{per}(\Gamma_{\V^*})$. 

Therefore, summing up all the contributions, we search for $\tilde{u} = u^*+\delta u\in H^{1}_{h,\Gamma_\D}(\Omega_{\V^*})$ and $\delta \V \in \Theta$ such that
\begin{IEEEeqnarray*}{rCl}
	\IEEEyesnumber\label{eq:ProtoWeak} \IEEEyessubnumber*
	\int_{\Omega_{\V^*}} \mkern-12mu \nabla \tilde{u} \cdot \nabla \varphi \d\Omega &+&
	\int_{\Gamma_{\V^*}} (\partial_n u^* - g) \, \partial_n \varphi \; \delta\V \cdot n \d \Gamma  
	+ \int_{\Gamma_{\V^*}} \mkern-12mu \nabla_\Gamma u^* \cdot \nabla \varphi \; \delta\V \cdot n \d \Gamma  \IEEEnonumber  \\
	&-& \int_{\Gamma_{\V^*}} \mkern-12mu \K_\H \varphi \; \delta\V\cdot n \d\Gamma
	= \int_{\Omega_{\V^*}} \mkern-12mu f\,\varphi \d\Omega + \int_{\Gamma_{\V^*}} \mkern-12mu g\,\varphi \d\Gamma \label{eq:ProtoWeak1}
\end{IEEEeqnarray*}
\begin{IEEEeqnarray*}{rcl}
	\IEEEyessubnumber*
	\int_{\Gamma_{\V^*}} \mkern-12mu \tilde{u}\, v \d\Gamma +
	&\int_{\Gamma_{\V^*}} & \mkern-12mu \left[ \left(\partial_n u^* -\partial_n h + \H (u^* - h) \right) v + (u^* - h) \, \partial_n v \right]\, \delta \V \cdot n \d\Gamma \IEEEnonumber \\
	&=&  \int_{\Gamma_{\V^*}} \mkern-12mu h\, v \d\Gamma\,   \label{eq:ProtoWeak2}
\end{IEEEeqnarray*}
for all $\varphi \in H^1_{0,\Gamma_{\D}}(\Omega_{\V^*})$ and $v \in H^1_{per}(\Gamma_{\V^*})$.

In the above steps we used the \emph{tangential gradient splitting}, with the tangential gradient of a real function being defined as $\nabla_\Gamma(\cdot) = \nabla (\cdot) - \partial_n(\cdot) n $.

So far we carried out the computations in full generality, and \eqref{eq:ProtoWeak} is an exact Newton scheme.
We now proceed to comment on, and apply, some simplifications.
\begin{simplif}
	\label{simpl:constData}
	Without loss of generality one can consider $\partial_n h = 0$ on $\Gamma_{\V^*}$.
	Furthermore, we consider the case of constant data $h = h_0$, so then $\nabla_\Gamma h=0$ and $\nabla h = 0$ on $\Gamma_{\V^*}$.
\end{simplif}

\begin{simplif}
	\label{simpl:solution}
	The above formulas can be simplified further by considering, on $\Gamma_{\V^*}$, $u^*=h_0$ and $\partial_n u^* = g$.
	These conditions are consistent with the \emph{exact solution} of the FBP, and lead to a quasi-Newton method  as in \cite{Karkkainen1999fss, Vanderzee2013snm}.
\end{simplif}

Applying the above simplifications produces the following quasi-Newton scheme: Search for $\tilde{u} \in H^{1}_{h,\Gamma_\D}(\Omega_{\V^*})$ and $\delta \V \in \Theta$ such that
\begin{IEEEeqnarray*}{rCl}
	\IEEEyesnumber\label{eq:weak} \IEEEyessubnumber*
	\int_{\Omega_{\V^*}} \mkern-12mu \nabla \tilde{u} \cdot \nabla \varphi \d\Omega - \int_{\Gamma_{\V^*}} \mkern-12mu \K_\H\, \varphi\, \delta\V\cdot n \d\Gamma &=& \int_{\Omega_{\V^*}} \mkern-12mu f \varphi \d\Omega + \int_{\Gamma_{\V^*}} \mkern-12mu g\,\varphi \d\Gamma \label{eq:weak1} \IEEEeqnarraynumspace\\[5pt]
	\int_{\Gamma_{\V^*}} \mkern-12mu \tilde{u}\, v \d\Gamma + \int_{\Gamma_{\V^*}} \mkern-12mu g \, v \, \delta \V \cdot n \d\Gamma &=& \int_{\Gamma_{\V^*}} \mkern-12mu h_0\, v \d\Gamma\,  \label{eq:weak2}
\end{IEEEeqnarray*}
for all $(\varphi, v) \in H^1_{0,\Gamma_{\D}}(\Omega_{\V^*})\times  H^1_{per}(\Gamma_{\V^*})$.

\begin{remark}
	The Simplification \ref{simpl:solution} above is the reason why the scheme \eqref{eq:weak} is not and exact Newton scheme, but only quasi-Newton method: The derivatives are not calculated in the current approximation, but rather they are an approximation of the derivatives at the exact solution.
	This has the consequence that \eqref{eq:weak} does not achieve quadratic convergence, but only superlinear.
\end{remark}

\subsection{Linearisation of the very-weak formulation}
We now want to derive a linearisation for (\ref{eq:N}) at an arbitrary approximated solution $(u^*,\V^*)$, where as before $u^*\in H^{1}_{per}(\Omega_{\V^*})$ and $\V^*\in \Theta$.
The G\^{a}teaux derivative of $\N$ at $(u^*,\V^*)$ with respect to $u$ in the direction $\delta u$ is given by
\begin{equation}
	\langle \, \partial_u \N[u^*, \mathcal{V}^*; \varphi ] , \delta u \,\rangle = -\int_{\Omega_{\V^*}} \mkern-12mu \delta u\,\Delta \varphi \;\d\Omega. \label{eq:veryWeakShape_u}
\end{equation}
The linearisation with respect to the vector field is again performed using the Hadamard formulas (\ref{eq:hadamard}):
\begin{IEEEeqnarray}{rCl}
	\langle \, \partial_\V \N[u^*, \V^*; \varphi ] , \delta \V \,\rangle &=& \int_{\Gamma_{\V^*}} \mkern-12mu \nabla h \cdot \nabla \varphi \; \delta \V \cdot n \d\Gamma - \int_{\Gamma_{\V^*}} \mkern-12mu (u^* - h) \Delta \varphi \,  \delta\V \cdot n \d\Gamma  \IEEEnonumber\\
	&& - \int_{\Gamma_{\V^*}} \mkern-12mu \left[\K_\H \, \varphi + g\,\partial_n \varphi \right] \delta\V \cdot n \d\Gamma.
	\label{eq:veryWweakShape_V}
\end{IEEEeqnarray}
A Newton step at the point  $(u^*,\V^*)$ has then the following form: Search for $\delta u \in H^1_{0,\Gamma_\D}(\Omega_{\V^*})$ and $\delta \V \in \Theta$ such that  
\begin{equation}
\langle \, \partial_u \N[u^*, \mathcal{V}^*; \varphi ] , \delta u \,\rangle + \langle \, \partial_\V \N[u^*, \V^*; \varphi ] , \delta \V \,\rangle  = -\N(u^*, \V^*; \varphi)\, ,
\end{equation}
for all $\varphi\in H^2_{0,\Gamma_\D}(\Omega_{\V})$.

Summing the various terms we then search for $\tilde{u}=u^*+\delta u \in H^1_{h,\Gamma_{\D}}(\Omega_{\V^*})$ and $\delta\V\in\Theta$  such that \color{black}
\begin{IEEEeqnarray}{l}
	\label{eq:veryWeakNotSimplified}
	\int_{\Omega_{\V^*}} \mkern-10mu ( h-\tilde{u})\,\Delta \varphi \d\Omega - \int_{\Gamma_{\V^*}}  \mkern-12mu \left[\K_\H \varphi + g\,\partial_n \varphi + (u^* - h) \Delta \varphi\right] \,\delta\V\cdot n \d\Gamma \nonumber \\
	+ \int_{\Gamma_{\V^*}} \mkern-12mu \nabla h \cdot \nabla \varphi \, \delta \V \cdot n \d\Gamma
	= \int_{\Gamma_{\V^*}} \mkern-12mu g \varphi \d\Gamma + \int_{\Omega_{\V^*}} \mkern-12mu f \varphi \d\Omega
	- \int_{\Omega_{\V^*}} \mkern-12mu \nabla h \cdot \nabla \varphi \d\Omega,   \IEEEeqnarraynumspace
\end{IEEEeqnarray}
for all $\varphi\in H^2_{0,\Gamma_\D}(\Omega_{\V})$.

We proceed to apply Simplifications \ref{simpl:constData} and \ref{simpl:solution}, thereby obtaining the followings quasi-Newton scheme: Search for $\tilde{u}\in H^1_{h,\Gamma_{\D}}(\Omega_{\V^*})$ and $\delta\V\in\Theta$  such that
\begin{IEEEeqnarray}{rCl}
	\int_{\Omega_{\V^*}} \mkern-12mu ( h -\tilde{u})\,\Delta \varphi \d\Omega \, &-& \int_{\Gamma_{\V^*}}  \mkern-12mu \left[\K_\H \varphi + g\,\partial_n \varphi \right] \,\delta\V\cdot n \d\Gamma \IEEEnonumber \\[3pt]
	&=& \int_{\Gamma_{\V^*}} \mkern-12mu g \varphi \d\Gamma + \int_{\Omega_{\V^*}} \mkern-12mu f \varphi \d\Omega - \int_{\Omega_{\V^*}} \mkern-12mu \nabla h \cdot \nabla \varphi \d\Omega\, , \label{eq:veryWeak} \IEEEeqnarraynumspace
\end{IEEEeqnarray}
for all $\varphi\in H^2_{0,\Gamma_\D}(\Omega_{\V})$.

As we pointed out above, we cannot yet employ this formulation to produce a numerical scheme; we need to extract the strong form implied by \eqref{eq:veryWeak} and then write a new weak formulation.
Using standard variational arguments one can see that such strong form is:
\begin{IEEEeqnarray}{rCl's}
	\IEEEyesnumber\label{eq:newStrong} \IEEEyessubnumber*
	-\Delta\tilde{u} &=& f \quad & in $\Omega_{\V^*}$ \label{eq:newStrong_1}\\
	\partial\tilde{u}_n - \K_\H \; \delta\V\cdot n &=& g & on $\Gamma_{\V^*}$ \label{eq:newStrong_2}\\
	\tilde{u} &=& h & on $\Gamma_{\D}$ \label{eq:newStrong_3}\\
	g\,\delta\V\cdot n &=& h_0 - \tilde{u} & on $\Gamma_{\V^*}. $\label{eq:newStrong_4}
\end{IEEEeqnarray}
Thanks to the initial requirement on $g$ not vanishing, one can solve \eqref{eq:newStrong_4} for $\delta \V \cdot n$, obtaining the boundary update formula
\begin{equation}
	\label{eq:update}
	\delta \V \cdot n = \frac{h_0 - \tilde{u}}{g}.	
\end{equation}
Substituting in \eqref{eq:newStrong_2} and using \eqref{eq:newStrong_1}--\eqref{eq:newStrong_3} allows to write the new weak formulation: Search for $\tilde{u} \in H^1_{h, \Gamma_\D}(\Omega_{\V^*})$ such that
\begin{equation}
	\label{eq:newWeak_galerkin}
	\int_{\Omega_{\V^*}} \mkern-12mu \nabla\tilde{u} \cdot \nabla \varphi \d \Omega - \int_{\Gamma_{\V^*}} \left( \K_\H \, \frac{h_0 - \tilde{u}}{g} + g \right) \varphi \d \Gamma = \int_{\Omega_{\V^*}} f \varphi \; \d\Omega,
\end{equation}
for all $\varphi \in H^1_{0, \Gamma_\D}(\Omega_{\V^*})$.

\begin{remark}
	Solving Equation \eqref{eq:weak2} for $\delta \V \cdot n$ one obtains exactly Equation \eqref{eq:update}.
	Plugging then into Equation \eqref{eq:weak2} gives Equation \eqref{eq:newWeak_galerkin}.
	This shows that the two methods, the coupled system \eqref{eq:weak} and the formulation \eqref{eq:newWeak_galerkin} with boundary update as in \eqref{eq:update}, are variationally equivalent, so we can expect the behaviours of these two approaches to be very similar.
	On the other hand, even though they are equivalent in an infinite-dimensional setting, the difference in the way the vector field is handled (as a coupled projection in the former case, or a splitting method in the latter case) may be reflected in the performances at the discretised level.
	This will indeed be the case, as our numerical tests illustrate.
\end{remark}

The strong form \eqref{eq:newStrong} will also be used in the implementation of a collocation scheme, outlined in the next section.
In passing, we comment that in the case of non-constant Dirichlet data on the free boundary, from Equation \eqref{eq:veryWeakNotSimplified} one could split the gradient of $h$ in the third integral in its tangential and normal component, and apply the tangential Green's identity \cite[p. 367]{Delfour2011sag}. See also \cite{Vanderzee2013snm} for details.


\section{Numerical Schemes}
\label{sec:schemes}

In our numerical tests we used two Galerkin methods, one arising from \eqref{eq:weak} and one from  \eqref{eq:newWeak_galerkin}.
The main difference between them is that from the former one obtains a coupled method, while the latter yields a decoupled splitting method.
We implemented, moreover, a collocation method to solve the strong form \eqref{eq:newStrong}.

\subsection{B-splines based Isogeometric analysis}
This section presents the essentials of B-splines. For more details  we refer the interested reader to any of the specialised books on the subject, for instance \cite{Farin1990cas}.

A \textit{knot vector}  is a set of non-decreasing points $\Xi  =\lbrace \xi_1 \leq \ldots \leq \xi_{n+p+1}\ \rbrace$ with $\xi_i \in \mathbb{R}$ and $n$   the number of basis functions of degree $p$ to be built. 

A knot vector is said to be \textit{open} if its first and last knots have multiplicity $p+1$, and in this case it is customary to take $\xi_1 = 0$ and $\xi_{n+p+1}=1$.
The maximum multiplicity of each internal knot can never exceed $p$.
A knot vector is said to be \textit{uniform} if the knots are equispaced; in this case it is common to take $\xi_1 = -p\tau$ and $\xi_{n+p+1} =  p\tau$, with $\tau$ the distance between two consecutive knots. 

Univariate B-splines functions can be defined using the Cox-de Boor recursion formulas \cite{Deboor2001apg} as follows:

\indent
for $p=0$:
\begin{IEEEeqnarray*}{r"l}
	\hat{\psi}_{i,0}(\xi)= \begin{cases}1 &  \xi_{i}\leq \xi<\xi_{i+1}\\
		0 & \textrm{otherwise}
	\end{cases}
\end{IEEEeqnarray*}
\indent
for $p \geq 1$:
\begin{IEEEeqnarray*}{r"l}
	\hat{\psi}_{i,p}(\xi)=\! \begin{cases}\dfrac{\xi-\xi_{i}}{\xi_{i+p}-\xi_{i}}\hat{\psi}_{i,p-1}(\xi)+\dfrac{\xi_{i+p+1}-\xi}{\xi_{i+p+1}-\xi_{i+1}}\hat{\psi}_{i+1,p-1}(\xi)  & \xi_{i}\leq  \xi<\xi_{i+p+1} \\
		0  & \textrm{otherwise}
	\end{cases}
\end{IEEEeqnarray*}
where we adopt the convention $0/0=0$.
A B-spline basis function is therefore a piecewise polynomial in every knot span and at the knots it achieves regularity $C^{p-l}$ where $l$ is the multiplicity of the knot.
We will always use internal knots of multiplicity one, in order to have maximal regularity.

\begin{figure}[tbp]
	\begin{center}
		\subfloat[]{\includegraphics[width = .45\textwidth]{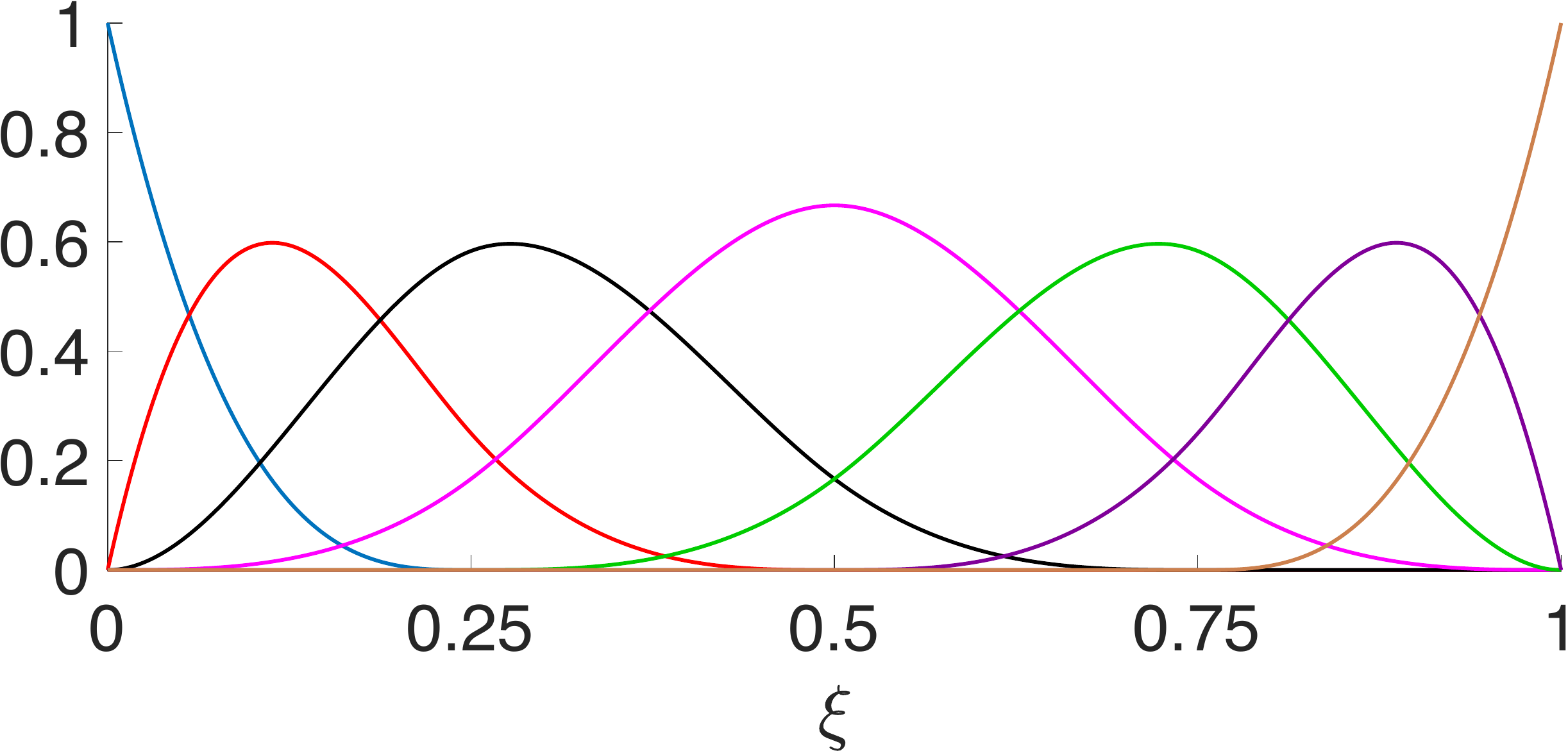}\label{fig:B-splines_open}}\quad
		\subfloat[]{\includegraphics[width = .45\textwidth]{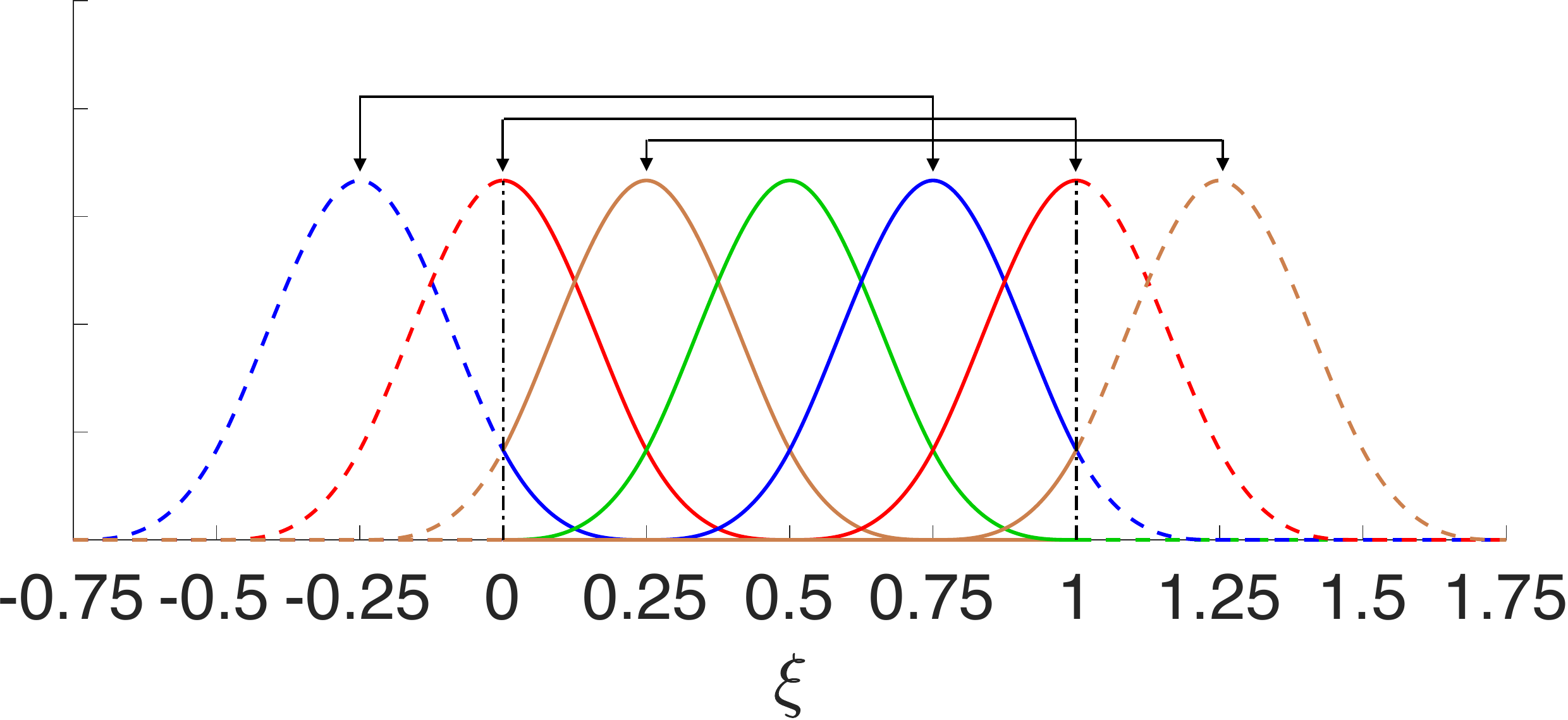}\label{fig:B-splines_periodic}} \qquad \qquad
	\end{center}
	\caption{Example of open and periodic B-spline basis. \textbf{(a)} Cubic basis on an open knot vector. \textbf{(b)} A periodic cubic basis on a uniform knot vector.}
\end{figure}

We denote with $\hat{S}^p= \text{span} \{\hat{\psi}_{i,p} \ | \ i=1, \ldots ,n \}$ the space spanned by $n$ B-splines of degree $p$.
We will often omit to explicitly indicate the polynomial degree.
On a uniform knot vector one can in addition construct a \textit{periodic} basis by appropriately identifying together functions laying at the beginning and at the end of the parametric domain:

\begin{equation}
\hat{S}_{per}^p= \mathrm{span} \{\hat{\psi}_{k}^{per} \} \quad \text{with } \quad
\begin{cases} \hat{\psi}_{k}^{per} := \hat{\psi}_{k} + \hat{\psi}_{n-p+k}, & k=1,\ldots,p ; \\
\hat{\psi}_{k}^{per} = \hat{\psi}_{k}, & \text{otherwise}
\end{cases}
\label{eq:periodic_space}
\end{equation}
Note that $\dim(\hat{S}_{per}^p)=n-p$.
Figure \ref{fig:B-splines_periodic} shows an example of maximum-regularity periodic B-splines basis with degree $p = 3$.

We can derive bivariate B-splines spaces, which we indicate in boldface, simply considering the tensor product of univariate ones.
Moreover, in our numerical tests we will use the same degree in each parametric direction.

Now, let $\mathbf{F}: \hat{\Omega}\rightarrow \Omega$ be a B-spline parametrisation of the physical domain $\Omega$,  and let $\hat{\mathbf{S}}^{p}$ be a space spanned by $N$ bivariate  B-splines $\hat{\phi}_k$ defined on the parametric domain $ \hat{\Omega}$.
Then, the corresponding space on $\Omega$ is defined as $\mathbf{S}^p=\mathrm{span}\{\phi_{k} \ | \ \phi_{k} = \hat{\phi}_k \circ \mathbf{F}^{-1}, \ k=1,\ldots,N\}$.
We moreover need to introduce a bivariate spline space spanned by functions periodic in $x$, that we denote $\mathbf{S}_{per}^p$.
This space is defined as the push-forward through the geometrical map $\mathbf{F}$ of the cross product between the periodic space $\hat{S}^p_{per}$, and the space $\hat{S}^p$ built from an open knot vector.

\subsection{Isogeometric Galerkin methods}

In both Galerkin-based schemes we choose as a trial space for $\tilde{u}$
\begin{equation}
{\mathbf{V}}_h^p:= {\mathbf{S}}^p_{per} \cap H^1_{h,\Gamma_\D}(\Omega_{\V^*}),
\label{eq:trial_space}
\end{equation} 
while as test space 
\begin{equation} 
{\mathbf{V}}_0^p:= {\mathbf{S}}^p_{per} \cap H^1_{0,\Gamma_\D}(\Omega_{\V^*}).
\label{eq:test_space}
\end{equation} 

The structure of the two algorithms is illustrated below. 

\begin{algorithm}
	\caption{ - Coupled Galerkin scheme}\label{al:coupledSplitting}
	\begin{algorithmic}[1]
		\State Choose the starting $\V_{0}$,
		\State Given $\V_k$, compute $(\tilde{u}_k, \delta\V\cdot n_k)$ solution of \eqref{eq:weak} in the domain $\Omega_{\V_{k}}$,
		\State Update the free boundary with $\V_{k+1} = \V_{k} + (\delta\V\cdot n_k)m_{k}$,
		\State Repeat steps 2--3 until $\norm{\delta\V\cdot n_k} \leq \mathrm{tol}$.
	\end{algorithmic}
\end{algorithm}

\begin{algorithm}
	\caption{ - Decoupled (splitting) Galerkin scheme}\label{al:decoupled}
	\begin{algorithmic}[1]
		\State Choose the starting $\V_{0}$,
		\State Given $\V_k$, compute $\tilde{u}_k$ solution of \eqref{eq:newWeak_galerkin} in the domain $\Omega_{\V_{k}}$,
		\State Compute $\delta\V\cdot n_k $ from \eqref{eq:update},
		\State Update the free boundary with $\V_{k+1} = \V_{k} + (\delta\V\cdot n_k) m_{k}$,
		\State Repeat steps 2--4  until $\norm{\delta\V\cdot n_k} \leq \mathrm{tol}$.
	\end{algorithmic}
\end{algorithm}

The vector field $m_k:\Gamma_{\V_k}\rightarrow\mathbb{R}$ represents the direction in which the update of the free boundary is performed, and has to satisfy $m_k\cdot n_k = 1$.
In our tests we choose to perform a vertical update, therefore selecting $m_k = \left[0,\, 1/( n_k )_y\right]$.
This choice allows to consider as unknown $\delta\V\cdot n$ instead of
$\delta\V$, which permits to discretise \eqref{eq:weak2} and
\eqref{eq:update} directly, using ${S}^p_{per}$ as both the test and trial space.
A choice of $m_k = n_k $ in the algorithms would instead amount to performing the update in the direction normal to the boundary.

\begin{remark}
	It is important to realise that when performing the update with Equation \eqref{eq:update} one has to divide two spline functions.
	The resulting function is therefore, in general, not a spline, and a projection onto the appropriate spline space is then required.
	In our tests we treated this by means of an $L^2$ projection into the space defined by the boundary test functions.
	After each boundary update, the internal mesh is then fitted using a Coons interpolation technique.
\end{remark}

\subsection{Isogeometric collocation method}
The isogeometric collocation method presented here is built from \eqref{eq:newStrong}: We solve \eqref{eq:newStrong_4} for $\delta\V\cdot n$ and  replace its value  in \eqref{eq:newStrong_2}, obtaining the following:

\begin{IEEEeqnarray*}{rCl"L}
	\IEEEyesnumber\label{eq:decoupledStrongForm} \IEEEyessubnumber*
	-\Delta\tilde u & = & f  &\mathrm{in}\  \Omega, \label{eq:decoupledStrongForm_1} \IEEEeqnarraynumspace\\
	\nabla\tilde{u}\cdot n-(\partial _n g + \text{H}\, g + f)\; \frac{h_0-\tilde{u}}{g} & = &  g & \mathrm{on} \ \Gamma_\V, \label{eq:decoupledStrongForm_2} \IEEEeqnarraynumspace\\
	\tilde{u} & = & h & \mathrm{on} \ \Gamma_\D, \label{eq:decoupledStrongForm_3} \IEEEeqnarraynumspace\\
	\delta\mathcal{V}\cdot n & = &  \frac{h_0-\tilde{u}}{g}  & \mathrm{on} \ \Gamma_{\V} \label{eq:decoupledStrongForm_4}.\IEEEeqnarraynumspace
\end{IEEEeqnarray*}
The structure of this algorithm is summarised below.

\begin{algorithm}
	\caption{ - Collocation scheme}\label{al:collocation}
	\begin{algorithmic}[1]
		\State Choose the starting $\V_{0}$,
		\State Given $\V_k$, compute $\tilde{u}_k$, collocated solution of \eqref{eq:decoupledStrongForm_1}--\eqref{eq:decoupledStrongForm_3},
		\State Compute $\delta\V\cdot n_k$ from \eqref{eq:decoupledStrongForm_4},
		\State Update the free boundary with $\V^{(k+1)} = \V^{(k)} +(\delta\V\cdot n_k)m_k$ ,
		\State Repeat steps 2--4 until $\norm{\delta\V\cdot n_k} \leq \mathrm{tol}$.
	\end{algorithmic}
\end{algorithm}

The solution of \eqref{eq:decoupledStrongForm_1}--\eqref{eq:decoupledStrongForm_3} and the boundary update \eqref{eq:decoupledStrongForm_4} are performed using a collocation approach.
Given the finite dimensional spaces ${\mathbf{V}}_{h}^p$  and ${S}^p_{per}$ in which we search for a solution $(\tilde{u},  \delta\V\cdot n)$, the idea is to accurately choose a number of points $\tau_1, \ldots, \tau_n \in \Omega$, called \emph{collocation points}, where $n$ is the number of degrees of freedom of the problem, and enforce the equations to hold strongly at those points.

\begin{figure}
	\centering
	\includegraphics[width = 0.8\textwidth]{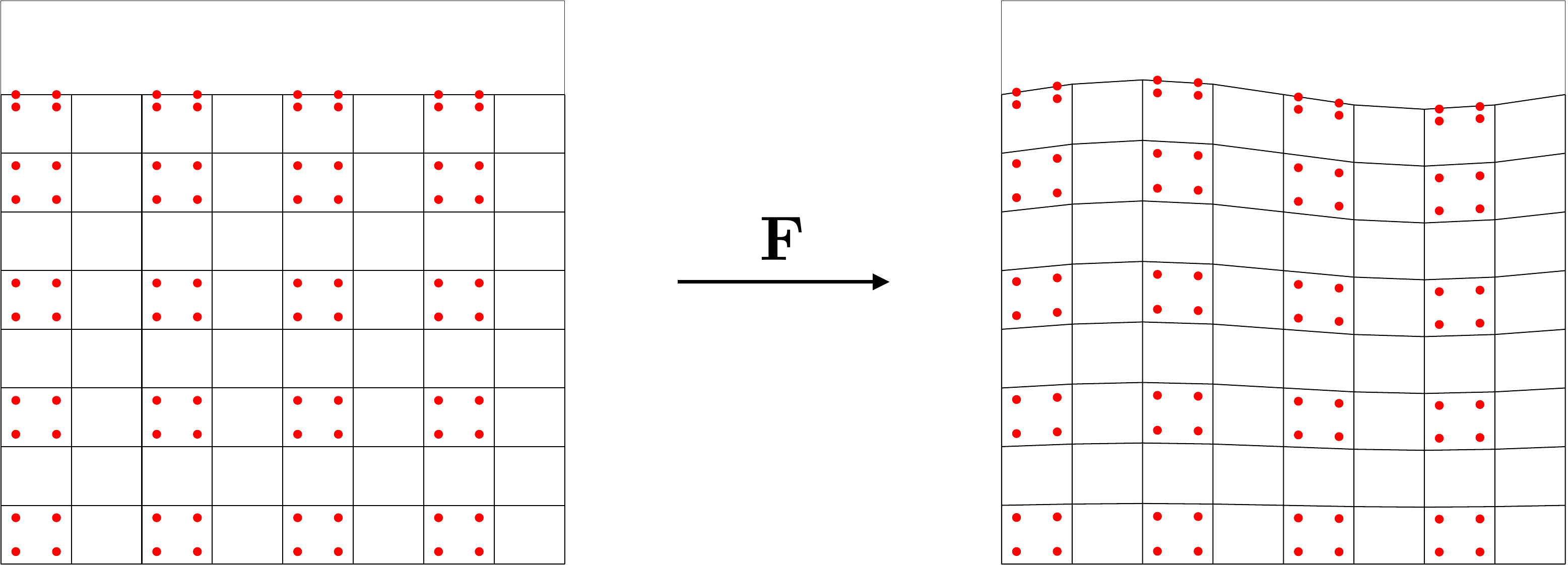}
	\caption{CSP collocation points in the parametric and in the physical domain. The points are the cross product of the periodic CSP points in the $x$-direction and the Dirichlet CSP points in the $y$-direction.}
	\label{fig:coll_pts}
\end{figure}

The appropriate selection of collocation points is crucial for the rate of convergence. 
Most of the classical choices of collocation points, for example, return suboptimal convergence rate even in a Poisson problem, contrary to the Galerkin approach which is optimal \cite{Schillingericc2013}.
However, the recent work \cite{Montardini2016ooi} suggests the use of a particular subset of Galerkin-superconvergent points, called clustered superconvergent points (CSP), as collocation points. 
This choice, that is the one that we adopt here, succeeds in achieving optimality for at least odd degrees B-splines discretisations.
In particular, the collocation points we use for the periodic problem \eqref{eq:decoupledStrongForm} are obtained by taking the cross product between univariate periodic CSP and univariate Dirichlet CSP (see \cite{Montardini2016ooi} for more details).
In our tests we however included also problems with only Dirichlet boundary conditions.
In that case the collocation points are selected as the push-forward of the cross-product of the univariate Dirichlet CSP points in the two parametric directions.
Figure \ref{fig:coll_pts} shows an example of CSP points in both the parametric and physical domain.
Note that we do not take any collocation points on the boundary $\{y = 0\}$, because we enforce the Dirichlet boundary conditions in the finite dimensional space that we consider, cf. \eqref{eq:trial_space}.

Similarly, the free boundary update is performed by collocating equation \eqref{eq:decoupledStrongForm_4} in the univariate periodic CSP, producing a fully-collocated scheme for problem \eqref{eq:decoupledStrongForm}.

\section{Numerical Results}
\label{sec:results}

This section collects our numerical results.
All algorithms have been implemented in Matlab using the GeoPDEs suite.
GeoPDEs is an Octave/Matlab software package for isogeometric analysis of partial differential equations \cite{Vazquez2016and}.
We applied the above Algorithms 1, 2, and 3 to different types of problems with either Dirichlet or periodic boundary conditions on the vertical sides.
It is clear that the error quantities in the problem are driven by the position of the free boundary: If the computed boundary matches the exact boundary solution, then the error on the internal function $u$ is simply the standard finite elements (IGA) or collocation approximation error.
For this reason, when evaluating the performance of the algorithms we have chosen the error quantities of interest to be the \emph{Dirichlet error}, $\norm{\tilde{u}(\Gamma_\V) - h_0}_{L^2}$, the error the computed function $u$ commits in satisfying the Dirichlet condition on the free boundary, and the \emph{surface position error}, $\norm{\Gamma_\V - \Gamma_{ex}}_{L^2}$, the error in the position of the computed free surface.

\subsection{Test 1: Parabolic boundary, Dirichlet b.c.}

This problem is constructed from the exact solution
\begin{equation}
\label{eq:exactSol}
u_{ex}(x,y) = \frac{y}{1+\alpha(x)} + \alpha(x)\, \frac{y}{1+\alpha(x)} \left( 1 - \frac{y}{1+\alpha(x)} \right)
\end{equation}
with
\begin{equation*}
\alpha(x) = \frac{1}{4}\, x\, (1-x).
\end{equation*}
The solution $u_{ex}$ attains constant value $u_{ex} |_{\Gamma_\V} = 1$ on the \emph{parabolic curve} $\Gamma_{ex} = \{ (x,y) \;|\; y = 1+\alpha(x), \, 0 \leqslant x \leqslant 1 \}$, which is therefore the exact free-boundary solution of the problem.

The data for problem \eqref{eq:FBP} are then found as follows:
\begin{IEEEeqnarray*}{rCl}
	f &=& -\Delta u_{ex}, \\
	g &=& \nabla u_{ex} \cdot \left(\textstyle\frac{1}{2}\,x - \frac{1}{4}, 1\right)/ \sqrt{1 + \textstyle\left(\frac{1}{2}\,x - \frac{1}{4}\right)^2}.
\end{IEEEeqnarray*}
We cast this problem with complete Dirichlet boundary conditions.
This amounts to imposing $h_0 = 1$ on the free boundary and $h=y$ on $\Gamma_\D \cup \Gamma_\P$.
We start our algorithms with $\Gamma_0 = \{ y = 1, \, 0 \leqslant x \leqslant 1 \}$ as an initial guess for the boundary.

\begin{figure}[htb]
	\centering
	\includegraphics[width=0.7\linewidth]{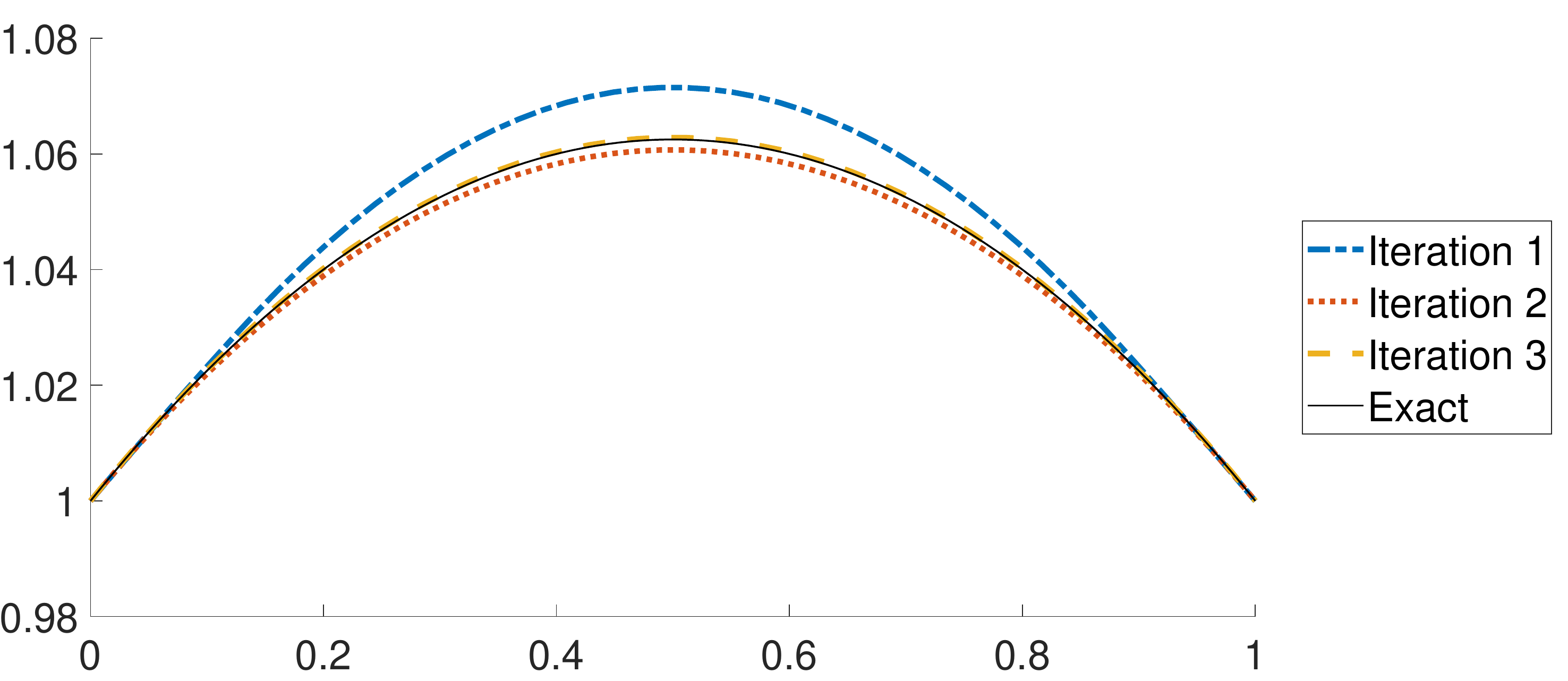}
	
	\caption{The first three iterations of Algorithm \ref{al:decoupled} for the Test 1 case, using a one element mesh and quadratic basis starting from a flat boundary with $y=1$.}
	\label{fig:parabolicBoundary}
\end{figure}

\begin{figure}[htb]
	\centering
	\subfloat{
		\includegraphics[height=2.8cm]{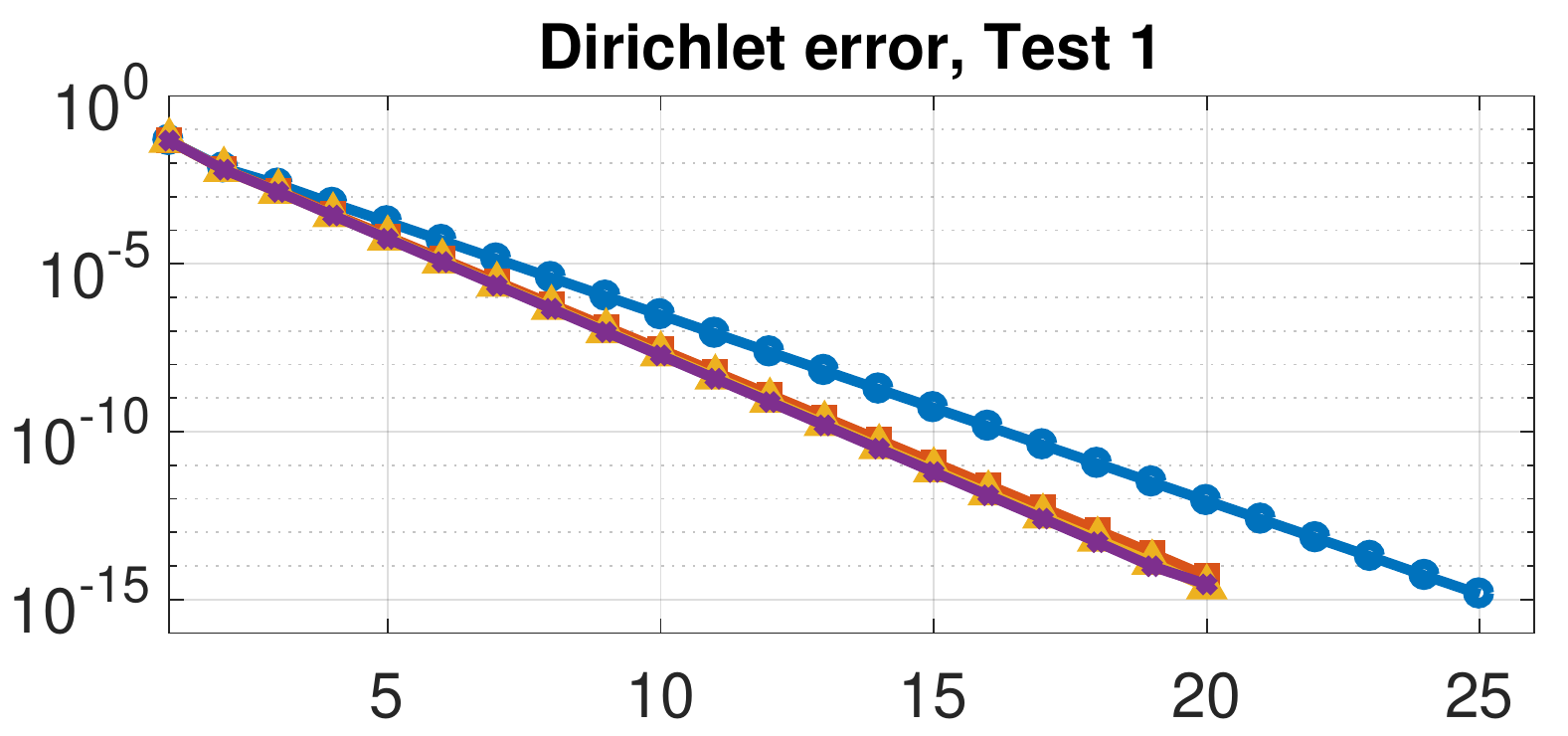}
	}
	\subfloat{
		\includegraphics[height=2.8cm]{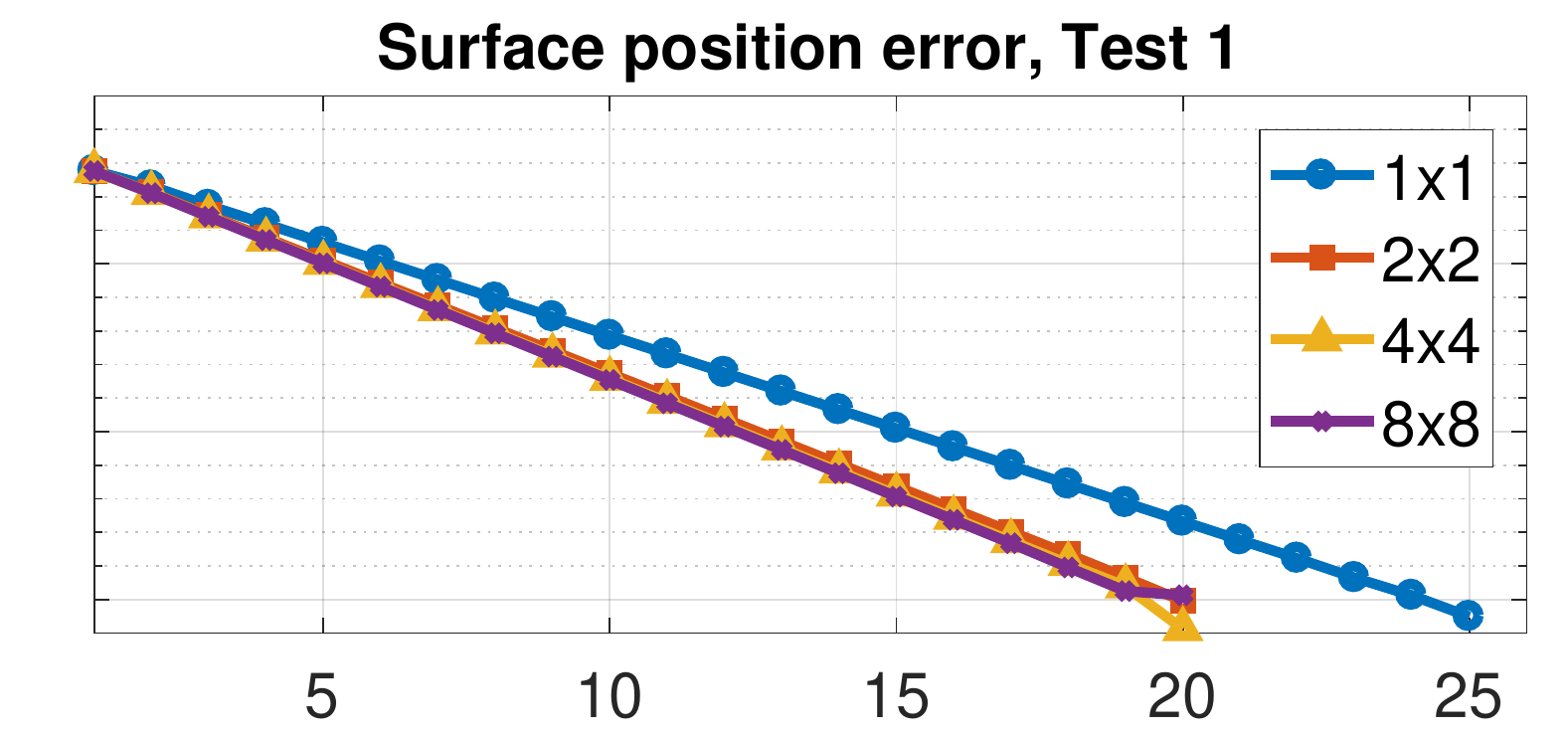}
	}
	
	\caption{Error quantities for Algorithm \ref{al:decoupled}, with a quadratic basis, on various mesh sizes. \textbf{(Left)} The Dirichlet error $\norm{\tilde{u}(\Gamma_\V) - h}_{L^2}$ as a function of the iterations. \textbf{(Right)} The surface position error $\norm{\Gamma_\V - \Gamma_{ex}}_{L^2}$. Machine precision is achieved for any mesh size. }
	\label{fig:errorTest1}
\end{figure}

Figure \ref{fig:parabolicBoundary} shows the first three iterations of the boundary update, together with the exact boundary solution, performed with a mesh with only 1 element and quadratic basis functions.
Those iterations have in particular been performed with Algorithm \ref{al:decoupled}, but Algorithm \ref{al:coupledSplitting} and Algorithm \ref{al:collocation} yielded identical results.
Figure \ref{fig:errorTest1} shows the convergence history of Algorithm 2 for both the Dirichlet error and the surface position error for various mesh sizes, using a quadratic basis.

\begin{figure}[htb]
	\centering
	\subfloat[Mesh 8x8]{
		\includegraphics[height=2.8cm]{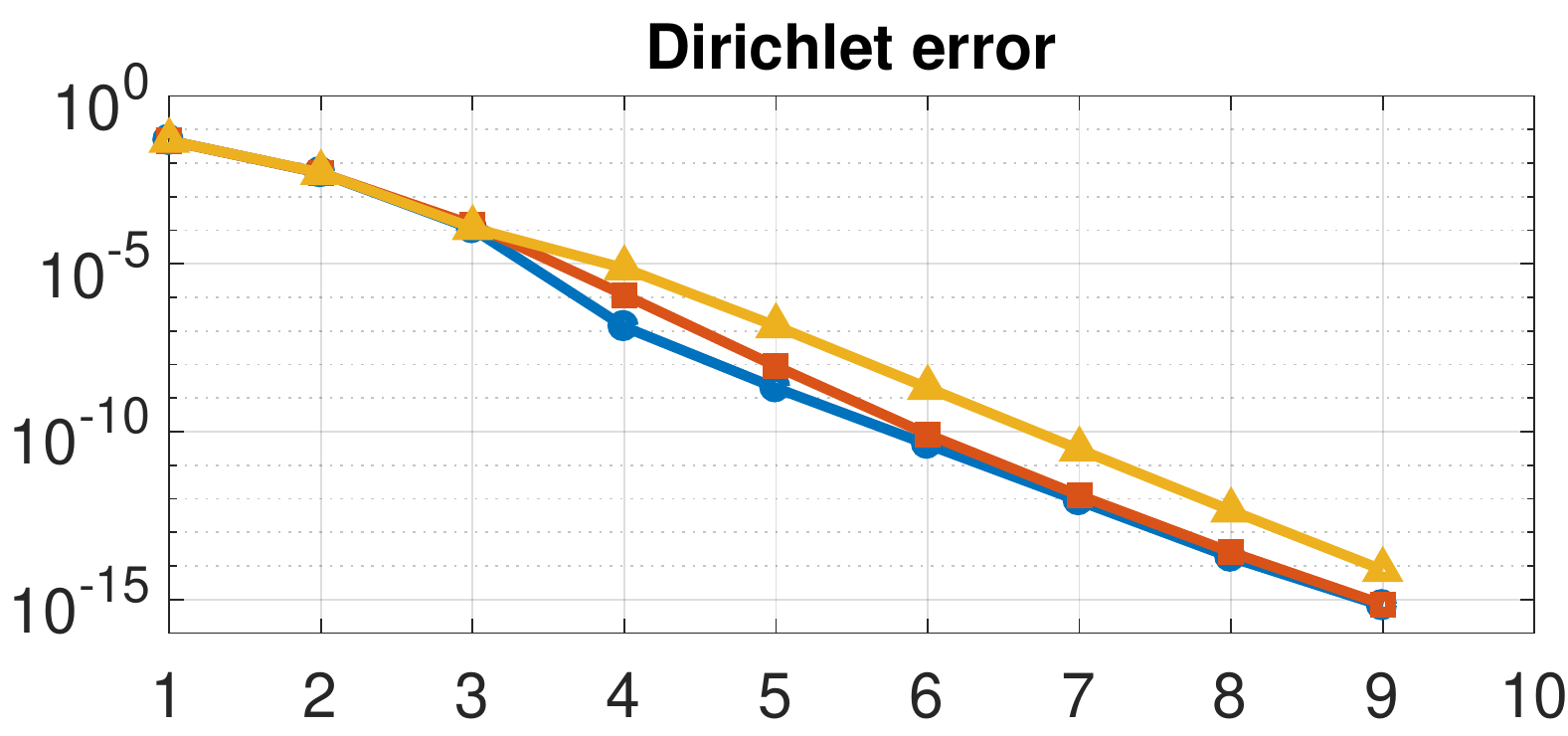}
		\hfill
		\includegraphics[height=2.8cm]{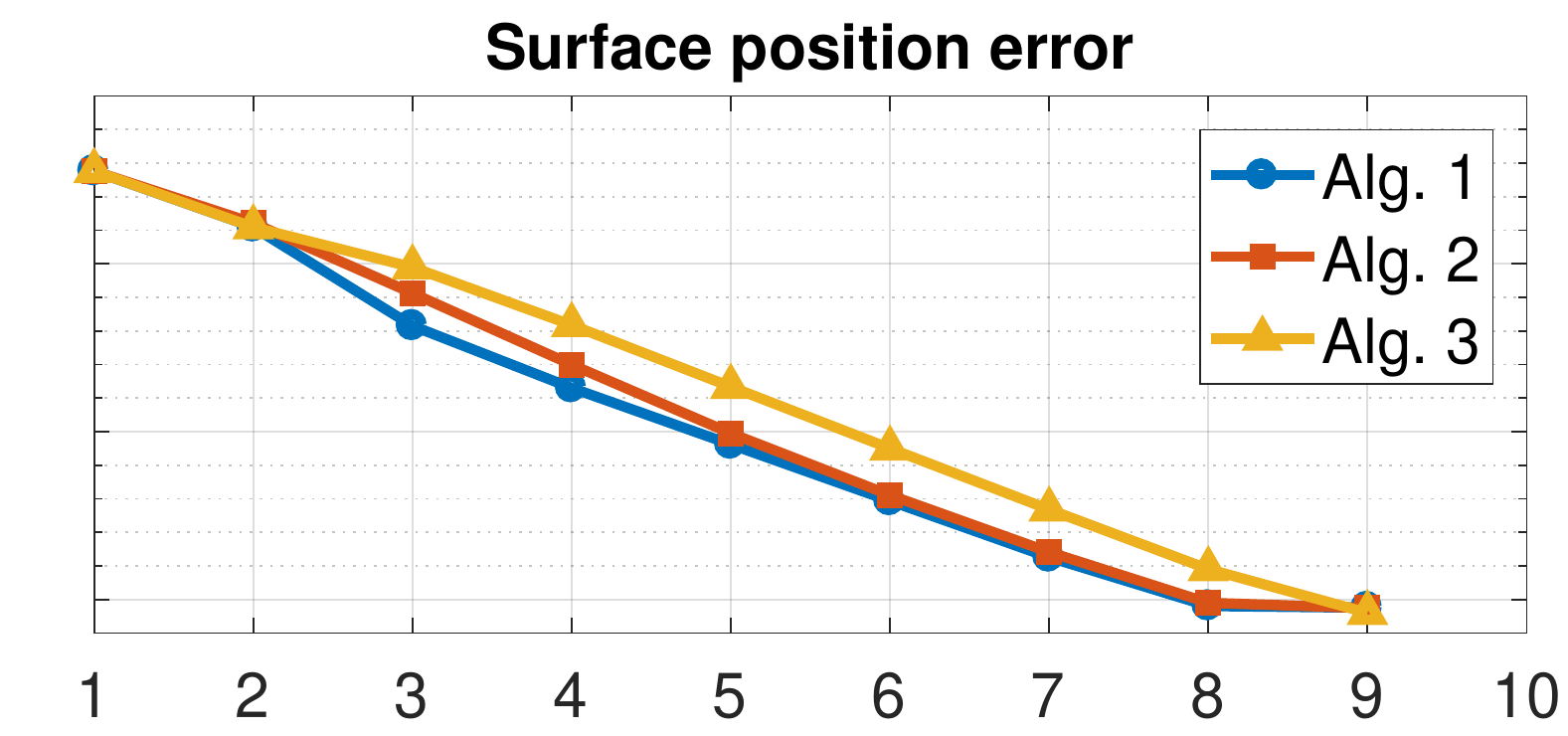}
	}
	
	\subfloat[Mesh 16x16]{
		\includegraphics[height=2.8cm]{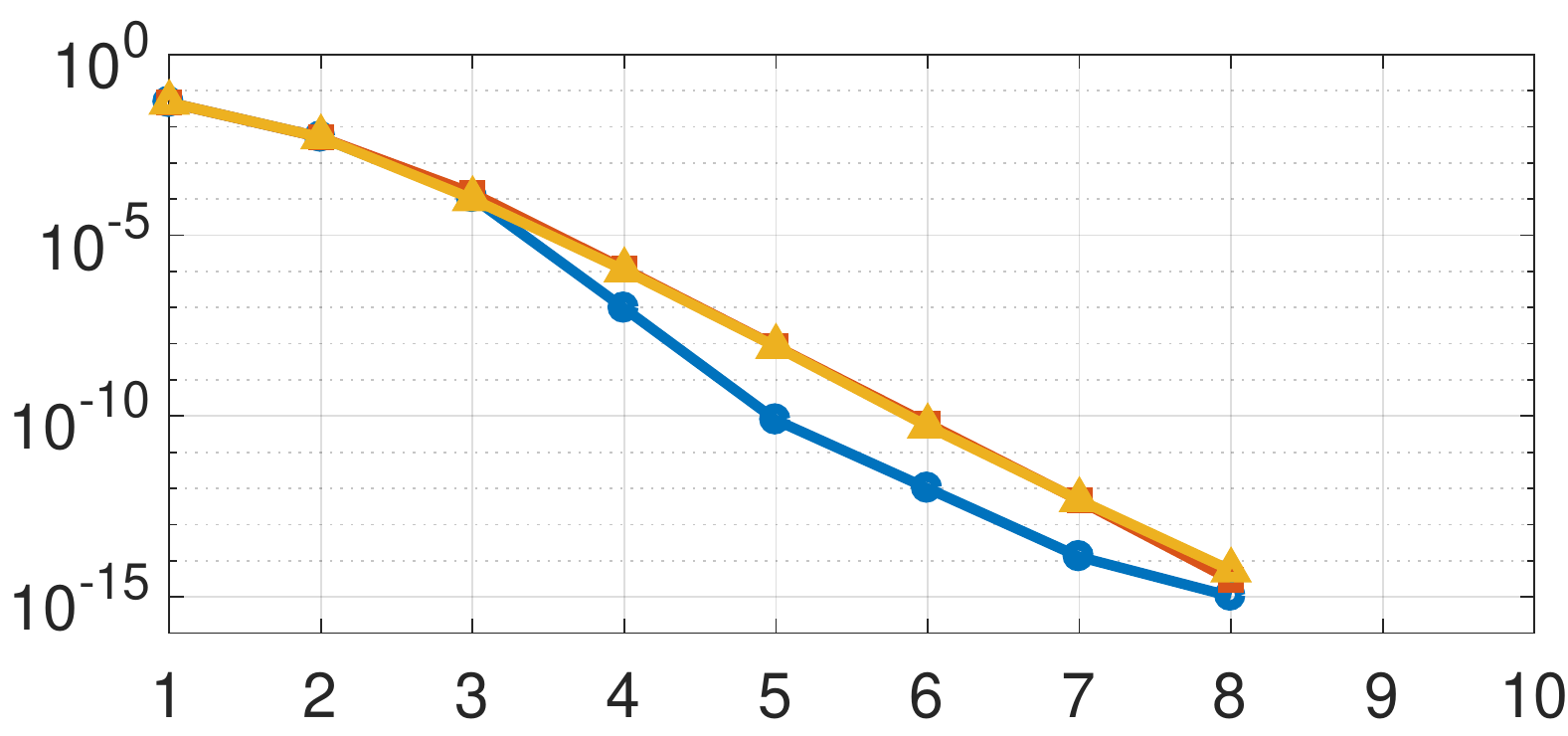}
		\hfill
		\includegraphics[height=2.8cm]{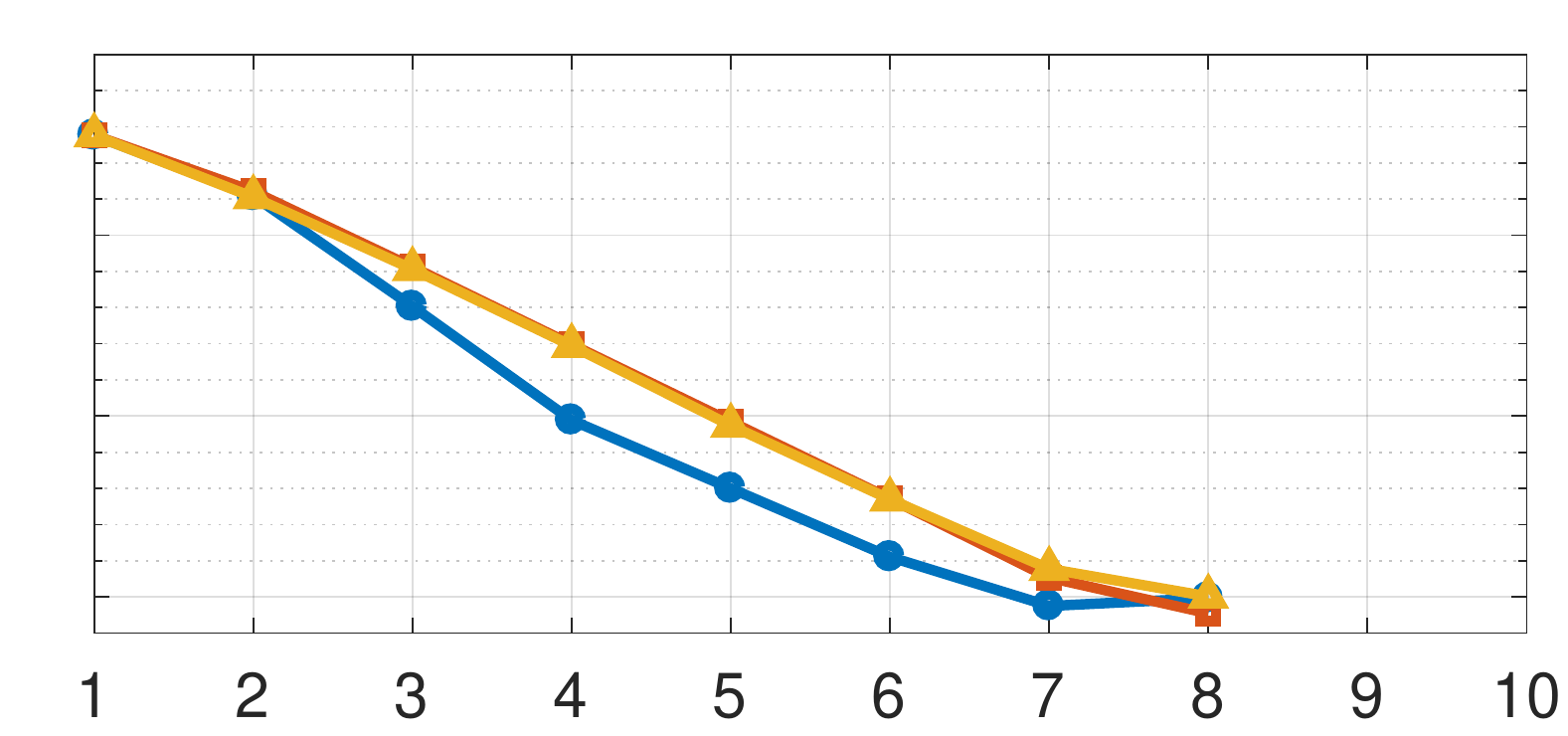}
	}
	
	\subfloat[Mesh 32x32]{
		\includegraphics[height=2.8cm]{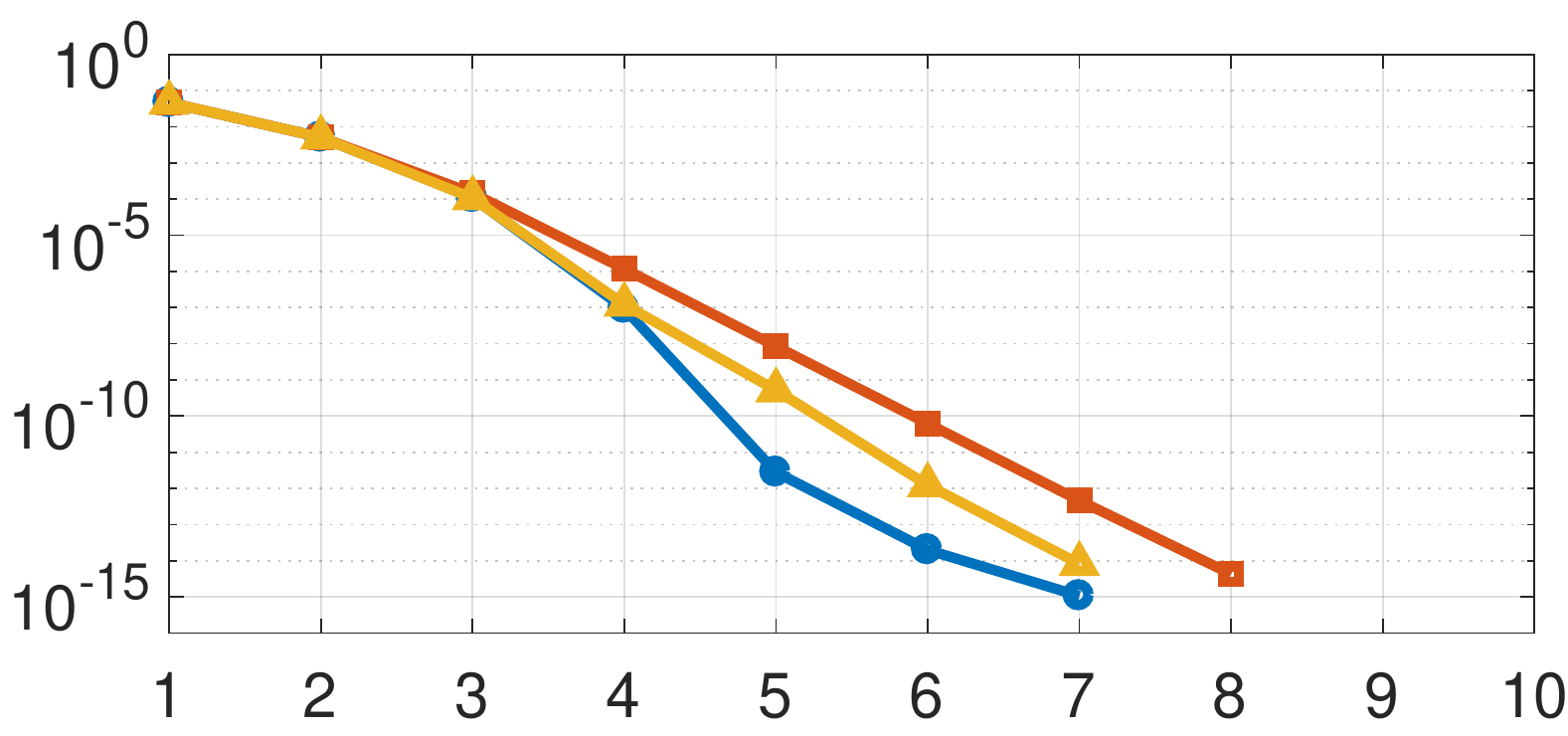}
		\hfill
		\includegraphics[height=2.8cm]{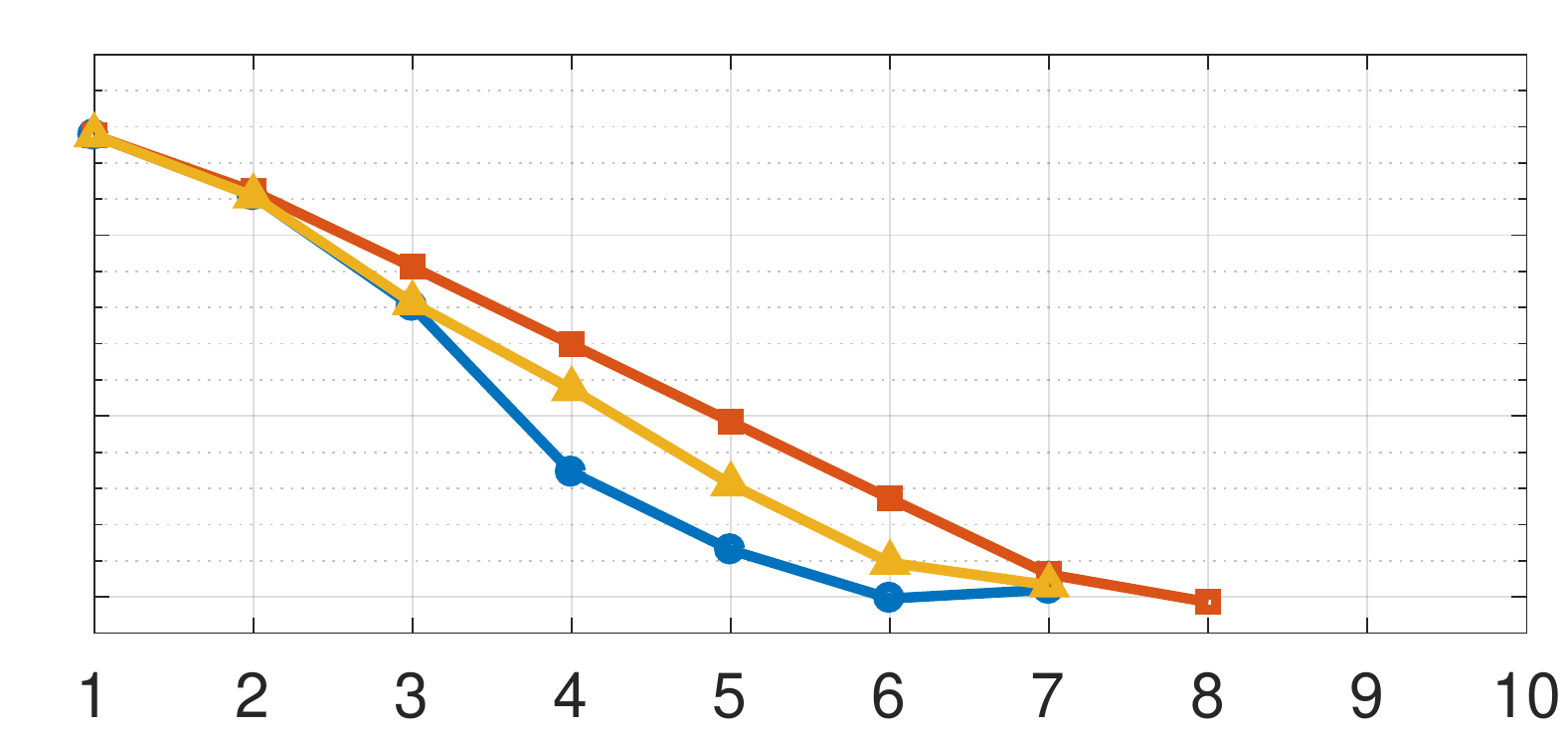}
	}
	
	\subfloat[Mesh 64x64]{
		\includegraphics[height=2.8cm]{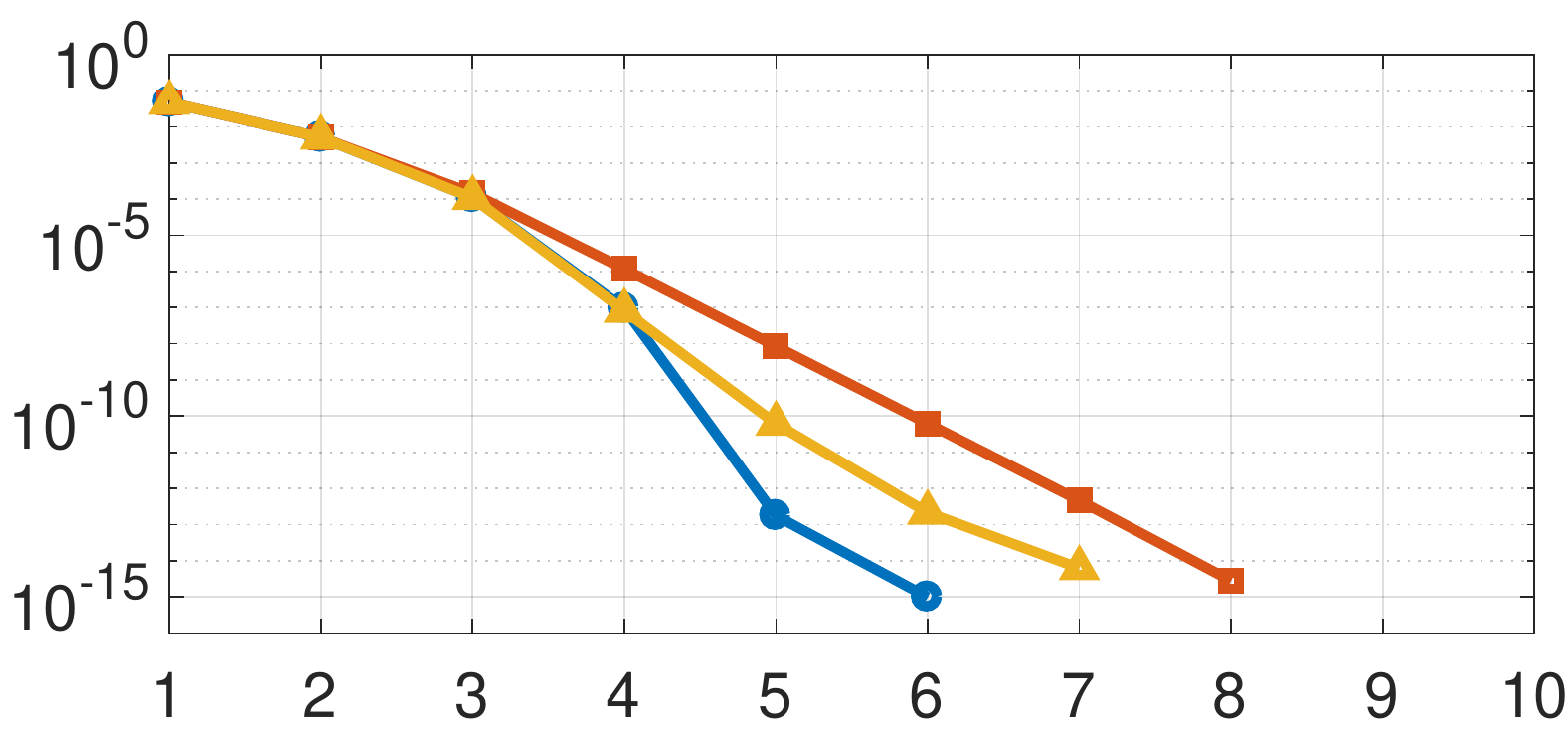}
		\hfill
		\includegraphics[height=2.8cm]{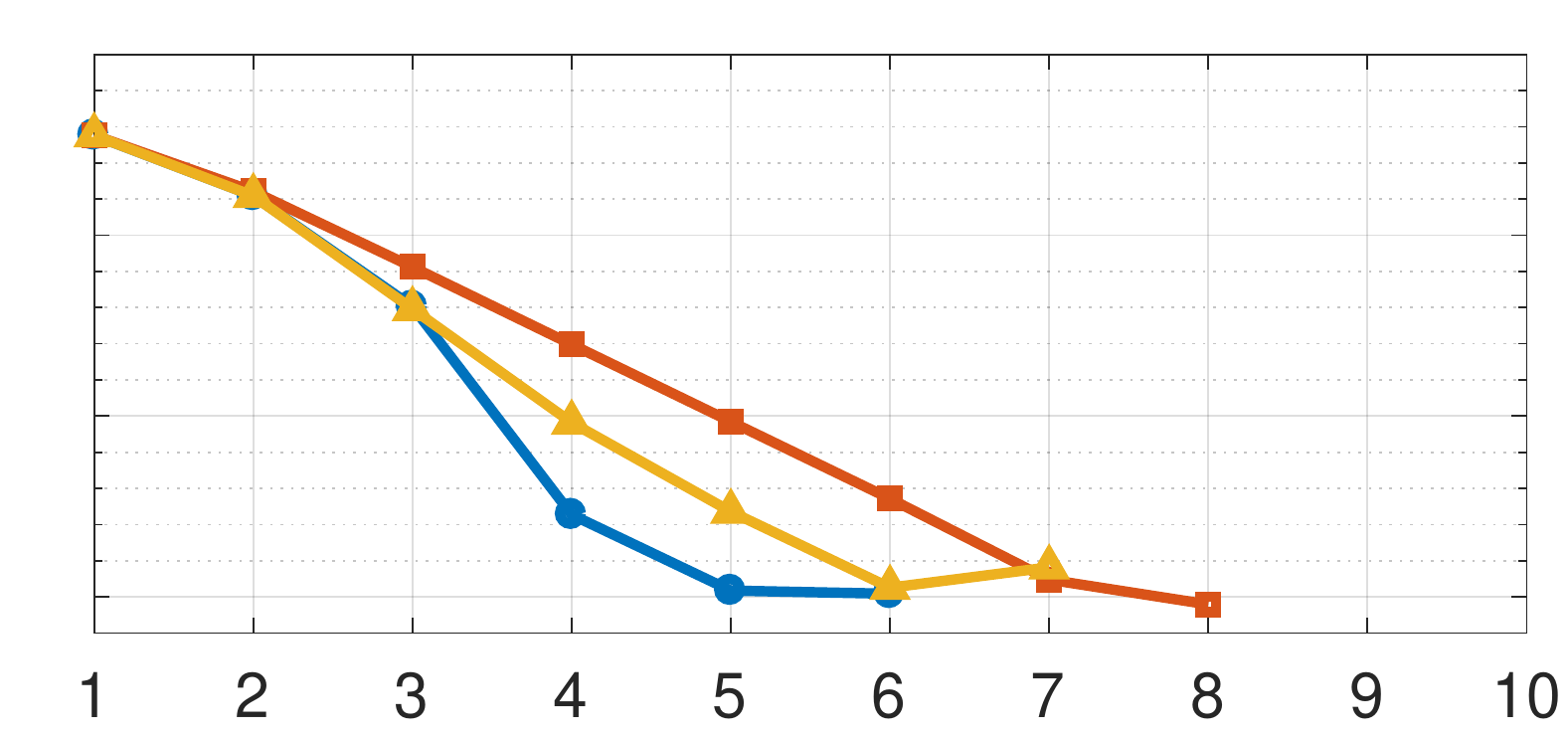}
	}

	\caption{A comparison of the three algorithms on Test 1 for different mesh sizes with cubic basis functions.}
	\label{fig:comparisonTest1}
\end{figure}

Figure \ref{fig:comparisonTest1} instead shows a comparison of the three different approaches using cubic basis functions.
The error plots show that Algorithm \ref{al:coupledSplitting}  improves the convergence speed once the solution is close enough.
The same behaviour is present also in the collocated scheme, Algorithm \ref{al:collocation}, albeit to a less degree, while it is not that apparent in Algorithm \ref{al:decoupled}.
However, all three algorithms' performances are quite similar on this test problem.
When it comes to runtime, Algorithm \ref{al:collocation} is much faster per iteration than the two Galerkin approaches, which is expected of a collocation scheme.

Note that this is the same setting as the ``Testcase I: Parabolic Free-Boundary'' presented in \cite[Section 5.2]{Vanderzee2013snm}.
However, in contrast to the results presented there we do not see a plateau in the error quantities, and machine precision is reached for any mesh size because the exact free boundary curve $\Gamma_{ex}$ and the exact solution
$u_{ex}$ restricted to $\Gamma_{ex}$ belong to the discrete space of the numerical approximation.

\FloatBarrier

\subsection{Test 2: Sinusoidal boundary, Dirichlet b.c.}

\begin{figure}[htb]
	\centering
	\includegraphics[width=0.7\linewidth]{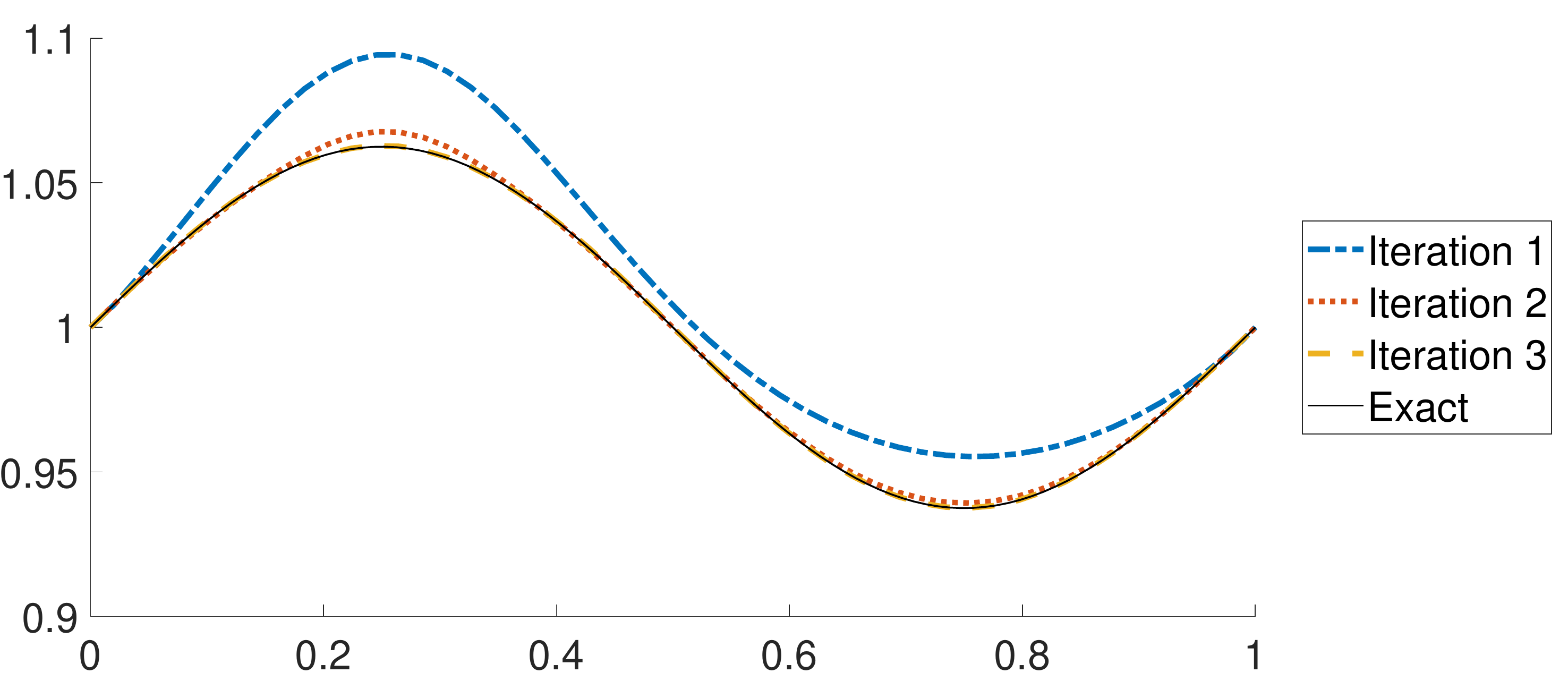}
	
	\caption{The first three iterations of Algorithm \ref{al:collocation} for the Test 2 case, with sinusoidal boundary and Dirichlet conditions, with an 8 elements mesh and cubic basis. Starting from a flat boundary with $y=1$.}
	\label{fig:sinusoidalBoundary}
\end{figure}

We now give an example where a plateau in the error is to be expected, and is actually found.
The problem data is derived as for Test 1 with an exact solution given by Equation \eqref{eq:exactSol} but with
\begin{equation*}
\alpha_{ex}(x) = \frac{1}{16}\, \sin(2\pi x),
\end{equation*}
so that the exact boundary $\Gamma_{ex} = \{ (x,y) \;|\; y = 1+\alpha(x), \, 0 \leqslant x \leqslant 1 \}$ is now a sinusoidal curve.
The boundary conditions are maintained of Dirichlet type, with $h_0 = 1$ on the free boundary, and $h=y$ on $\Gamma_\D \cup \Gamma_\P$.
Figure \ref{fig:sinusoidalBoundary} shows the first three boundary updates performed by Algorithm \ref{al:collocation}. The mesh is made of 8 elements, and the basis is cubic.
The initial boundary is again taken as the flat curve $\Gamma_0 = \{ y = 1, \, 0 \leqslant x \leqslant 1 \}$

\begin{figure}[htb]
	\centering
	\subfloat[Mesh 8x8]{
		\includegraphics[height=2.8cm]{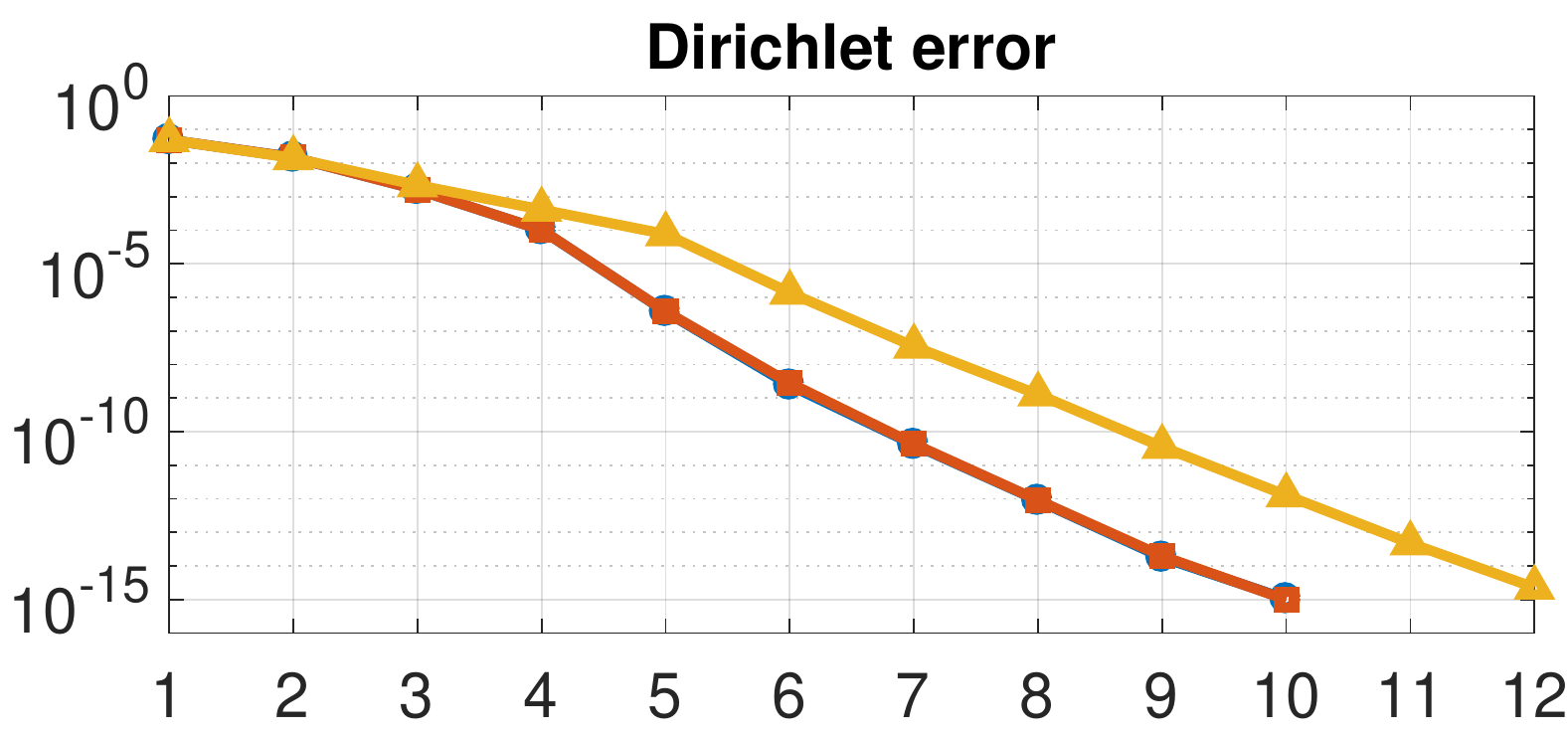}
		\hfill
		\includegraphics[height=2.8cm]{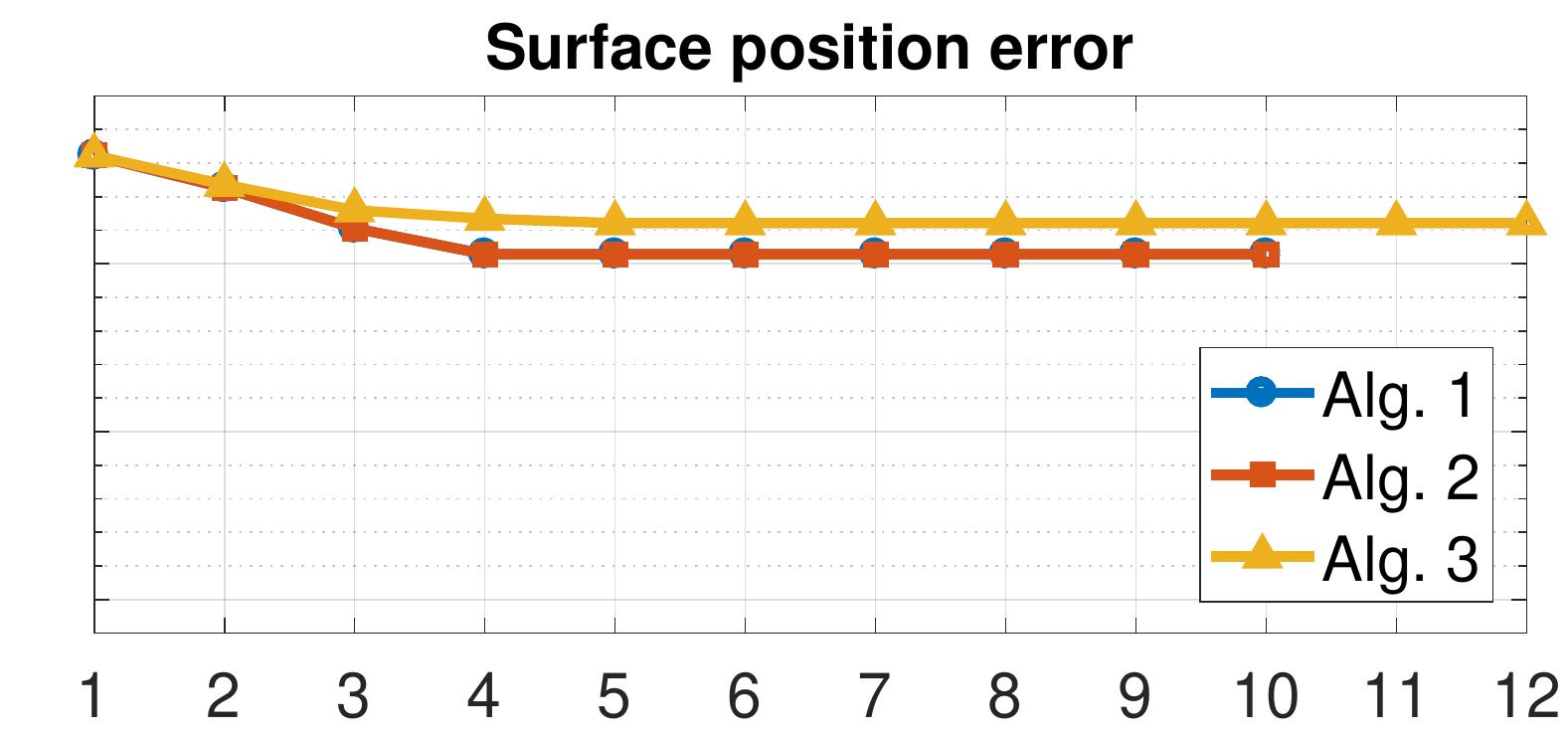}
	}
	
	\subfloat[Mesh 16x16]{
		\includegraphics[height=2.8cm]{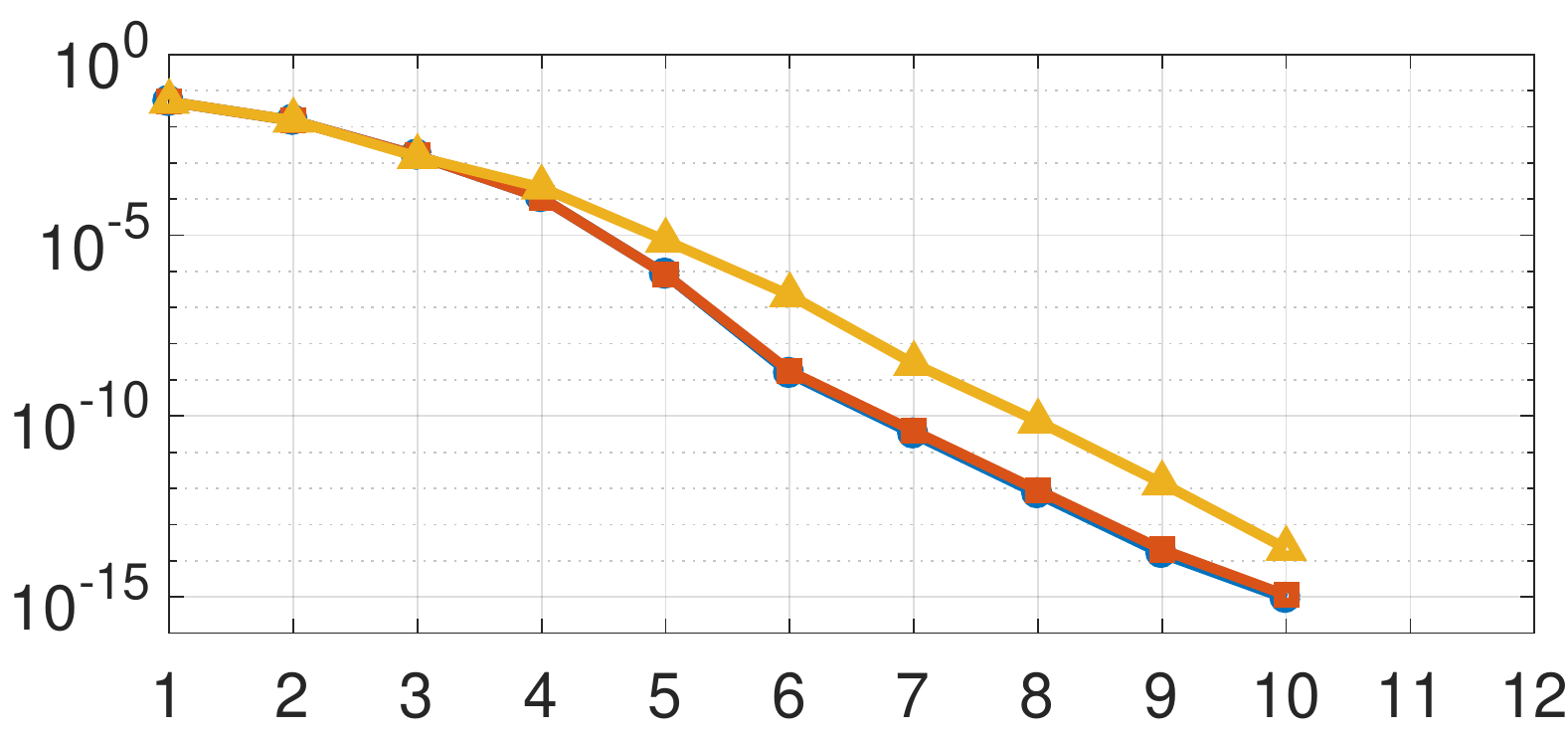}
		\hfill
		\includegraphics[height=2.8cm]{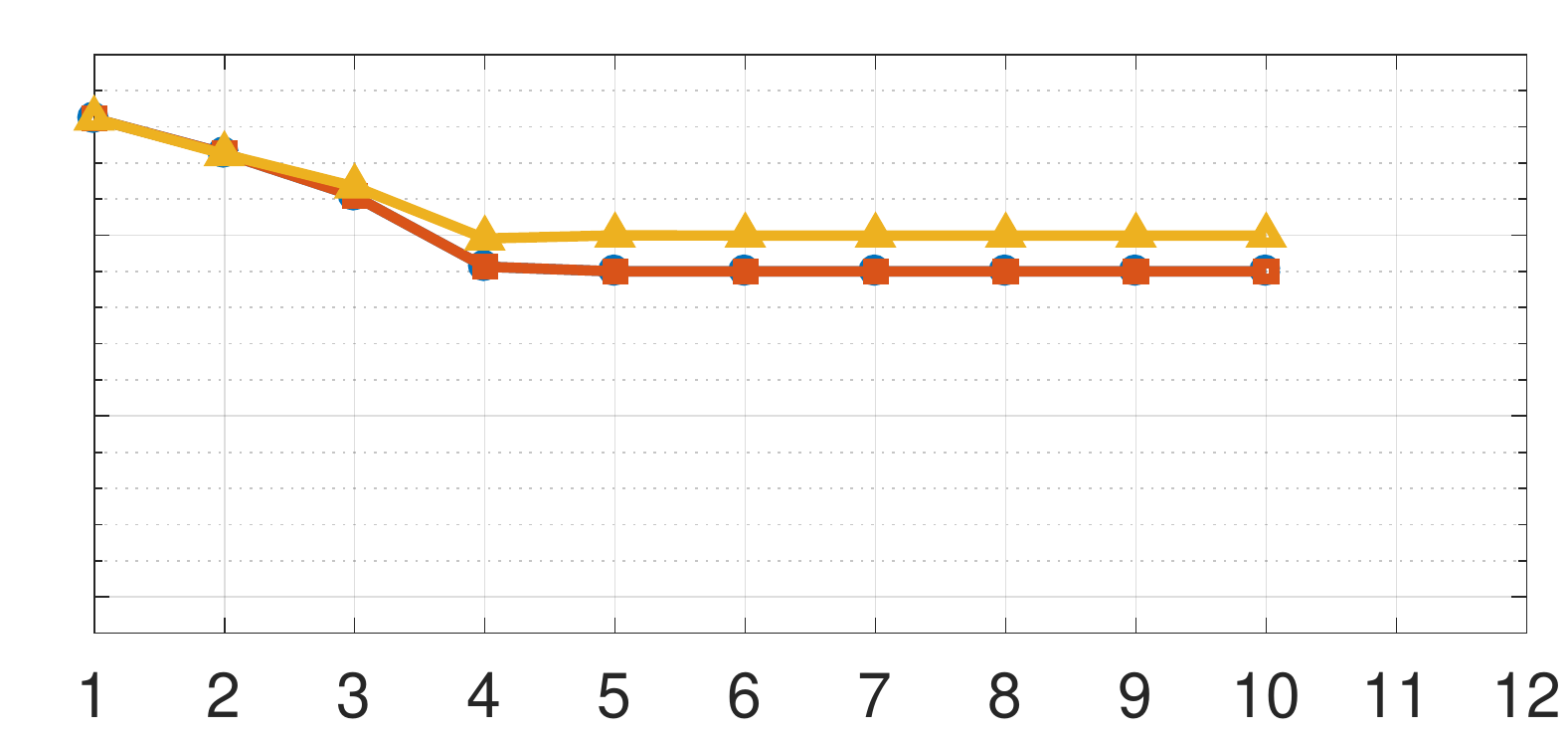}
	}
	
	\subfloat[Mesh 32x32]{
		\includegraphics[height=2.8cm]{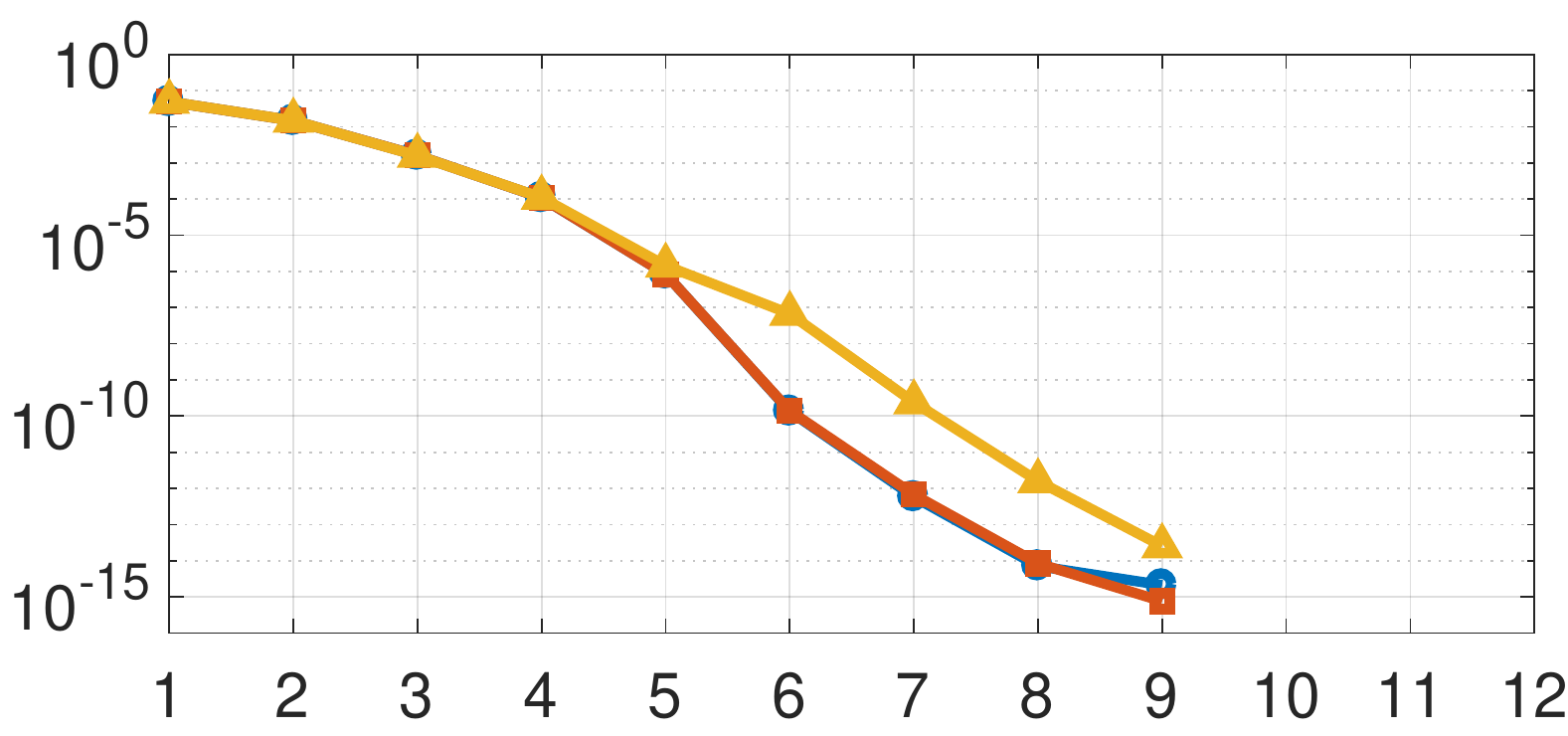}
		\hfill
		\includegraphics[height=2.8cm]{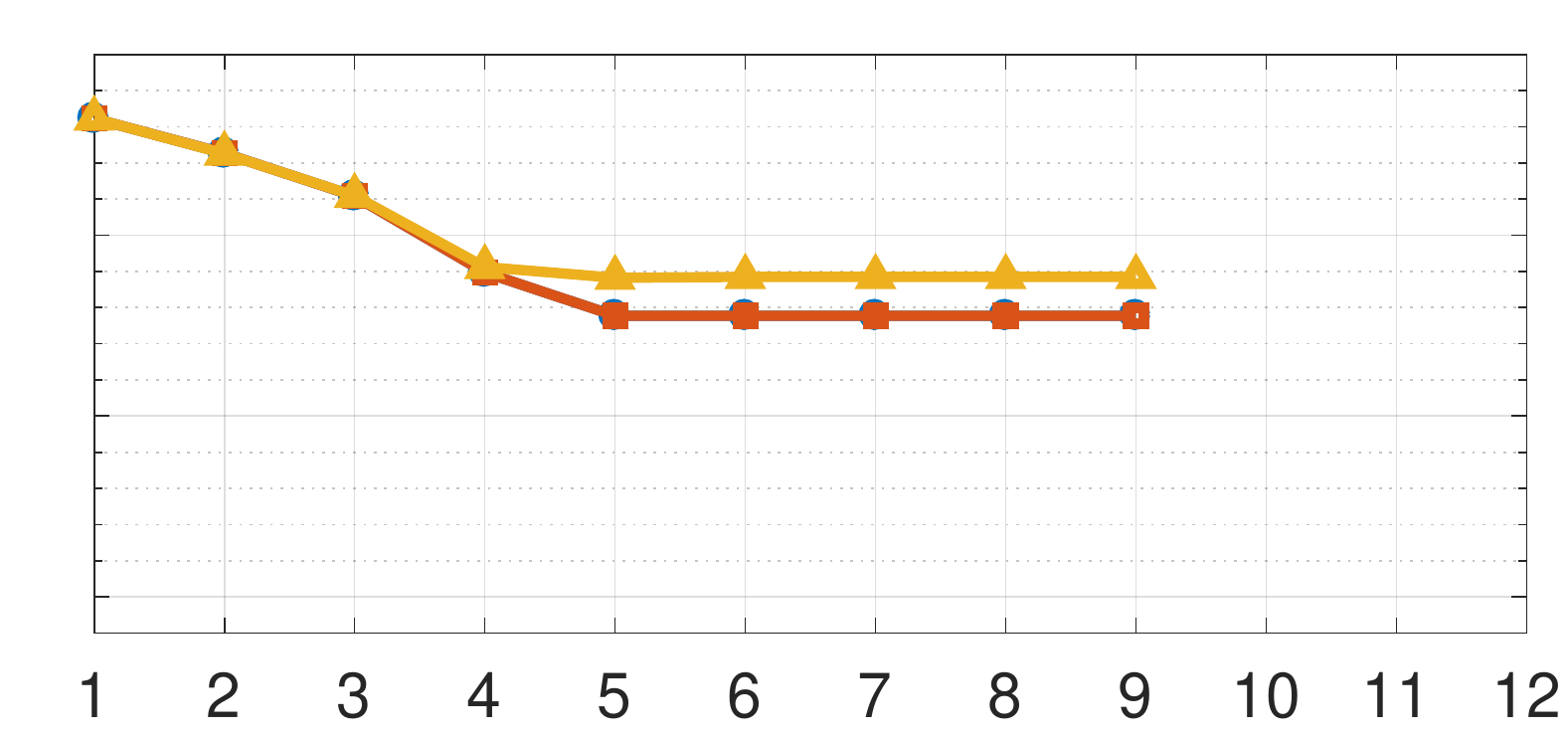}
	}
	
	\subfloat[Mesh 64x64]{
		\includegraphics[height=2.8cm]{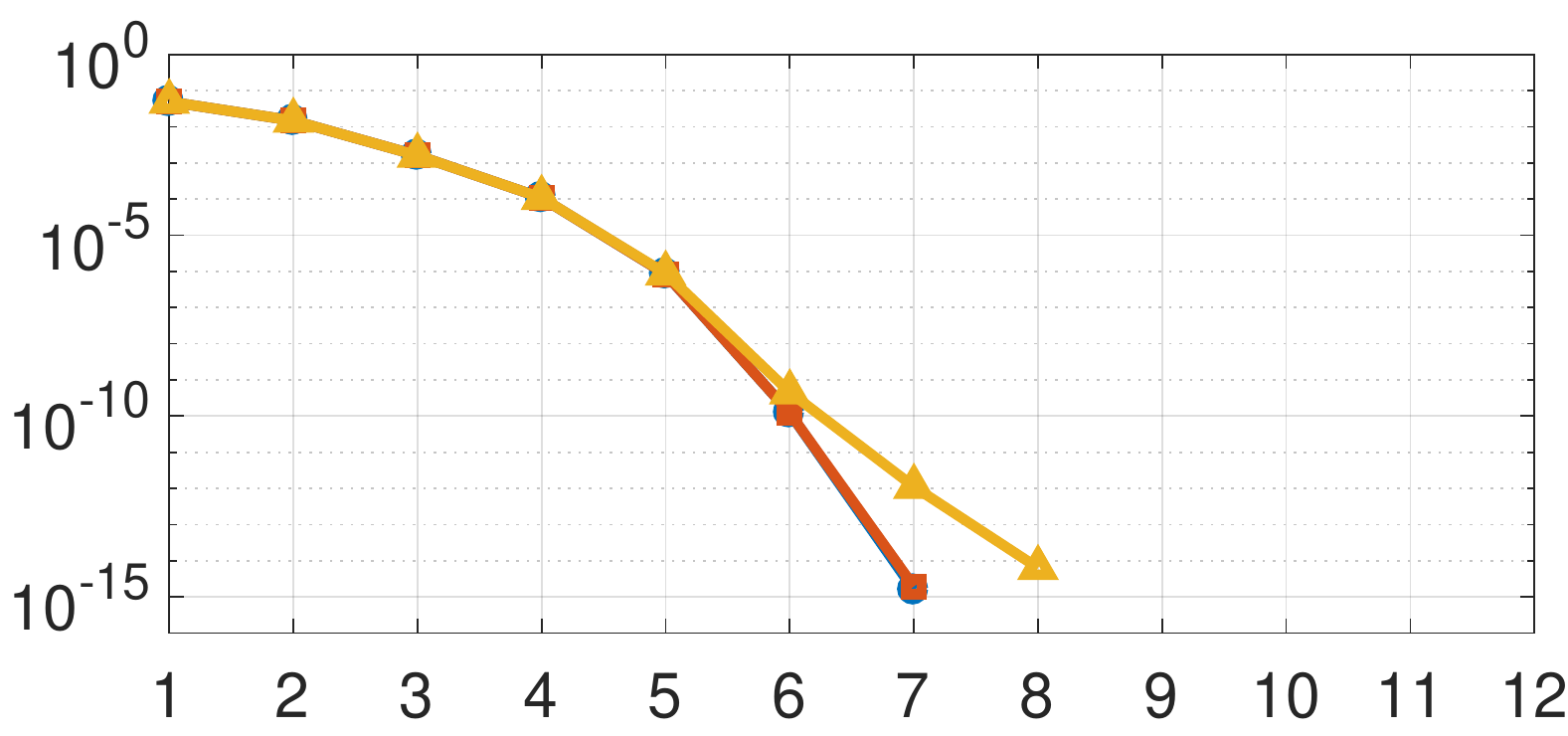}
		\hfill
		\includegraphics[height=2.8cm]{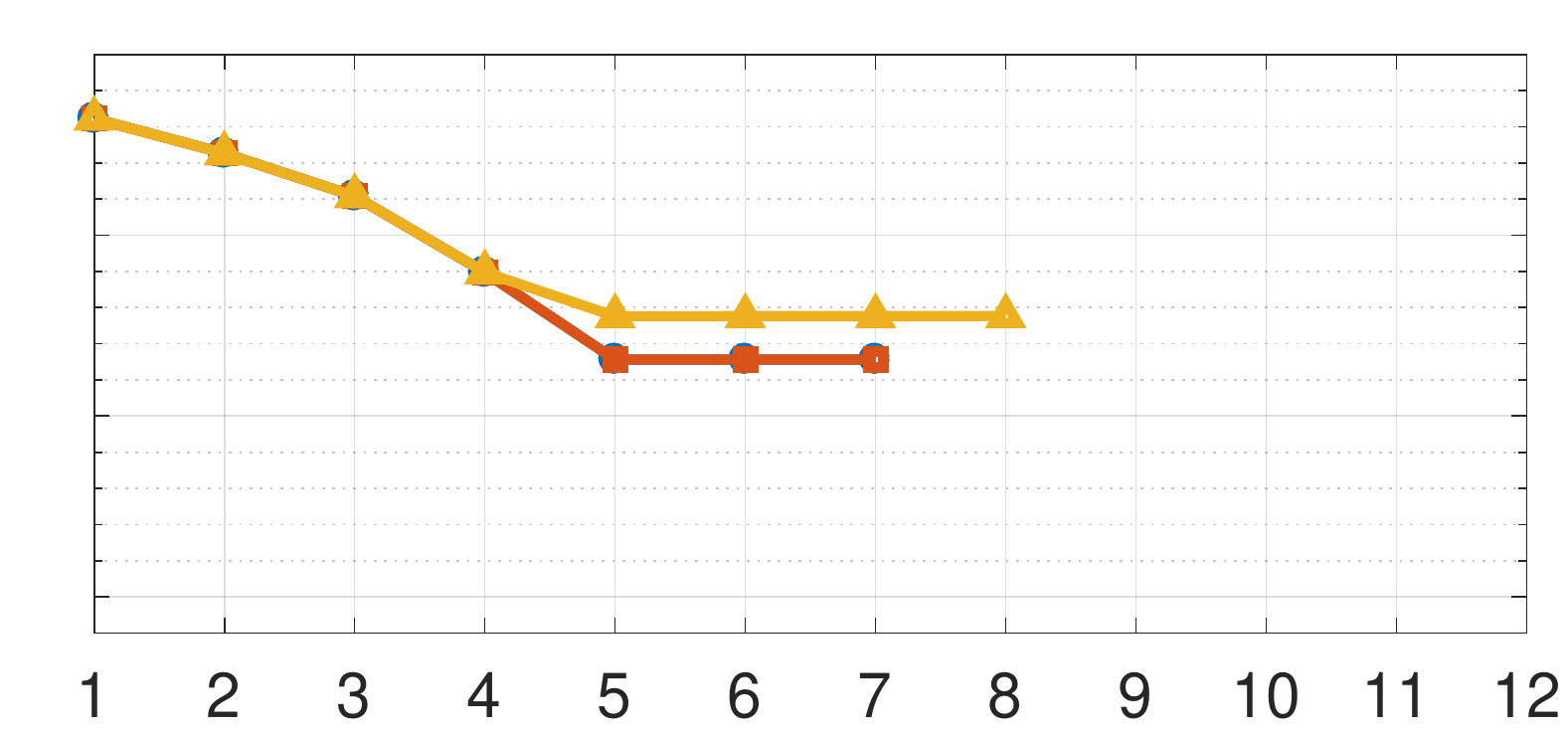}
	}

	\caption{A comparison of the three algorithms on Test 2 for different mesh sizes with cubic basis functions.}
	\label{fig:comparisonTest2}
\end{figure}

Figure \ref{fig:comparisonTest2} shows the error quantities vs iterations for the three algorithms.

As the mesh is refined we note that the collocation algorithm, Algorithm \ref{al:collocation}, has a slightly higher error than the other two approaches.
The surface position error, moreover, is abated with finer meshes in all approaches but remains always present.
This is due to the fact that a cubic B-spline cannot exactly represent a sinusoidal curve, and therefore the exact free boundary solution to this problems lies outside of the trial function space.
Lastly, Figure \ref{fig:comparisonTest2} shows how closely related Algorithms \ref{al:coupledSplitting} and \ref{al:decoupled} are, achieving almost identical performance on this benchmark test.

\FloatBarrier

\subsection{Test 3: Sinusoidal boundary, periodic b.c.}

In our third benchmark we employ the same problem data as in Test 2, but now periodic boundary conditions are placed on the lateral sides instead of Dirichlet ones.
In this test case we used the highest-possible regularity for the periodic conditions, meaning that the boundary functions are ``glued'' together with $C^{p-1}$ continuity.

\begin{figure}[htb]
	\centering
	\subfloat[Mesh 8x8]{
		\includegraphics[height=2.8cm]{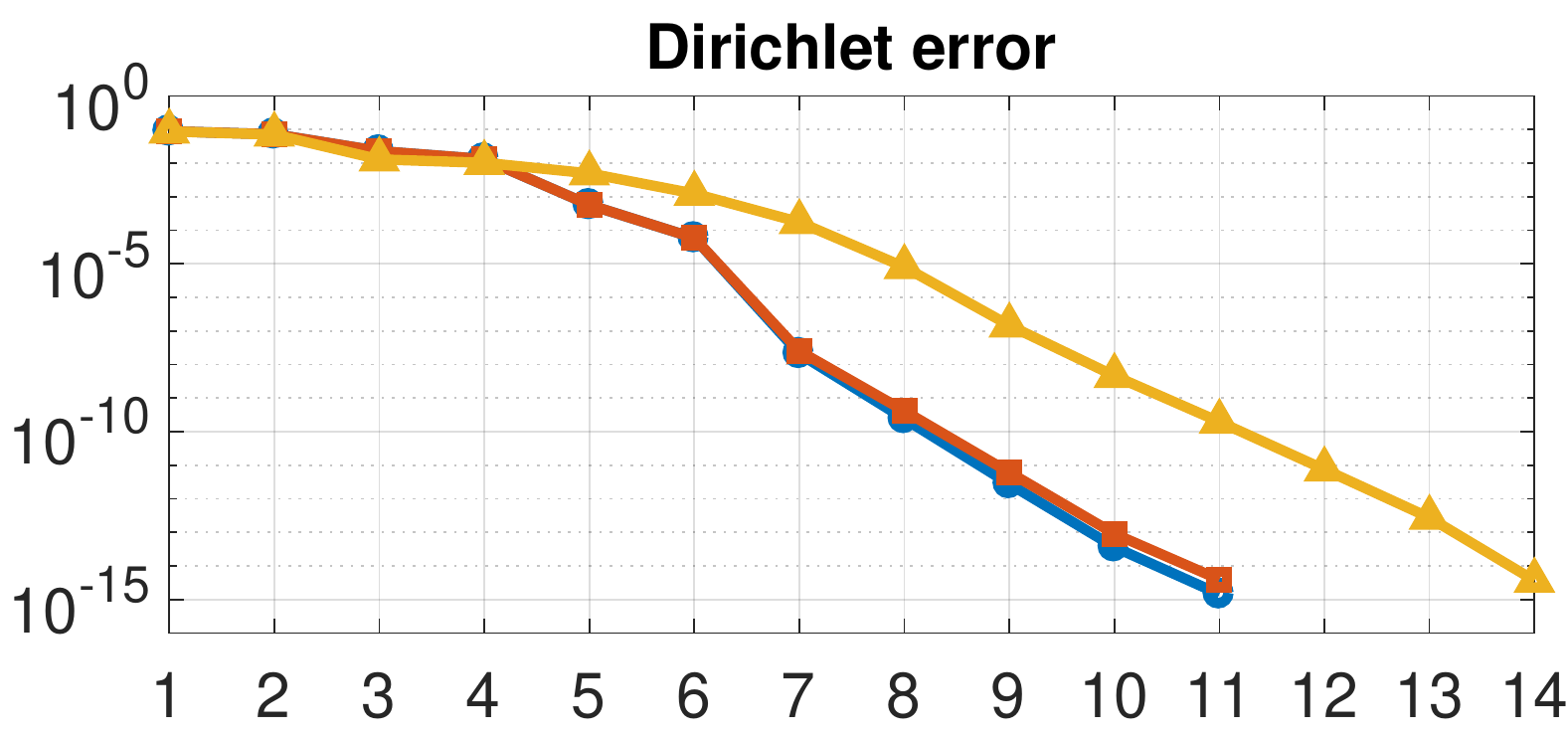}
		\hfill
		\includegraphics[height=2.8cm]{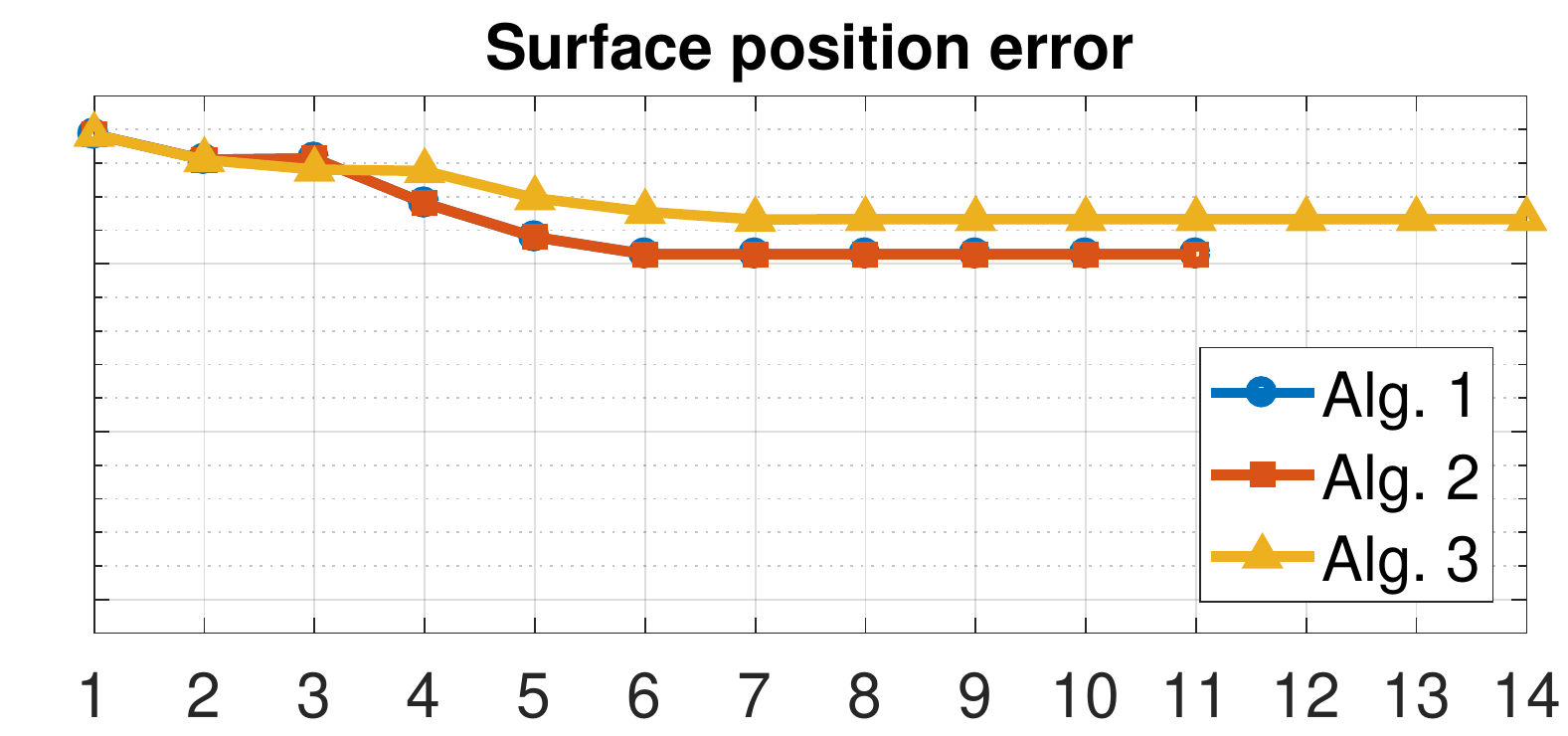}
	}
	
	\subfloat[Mesh 16x16]{
		\includegraphics[height=2.8cm]{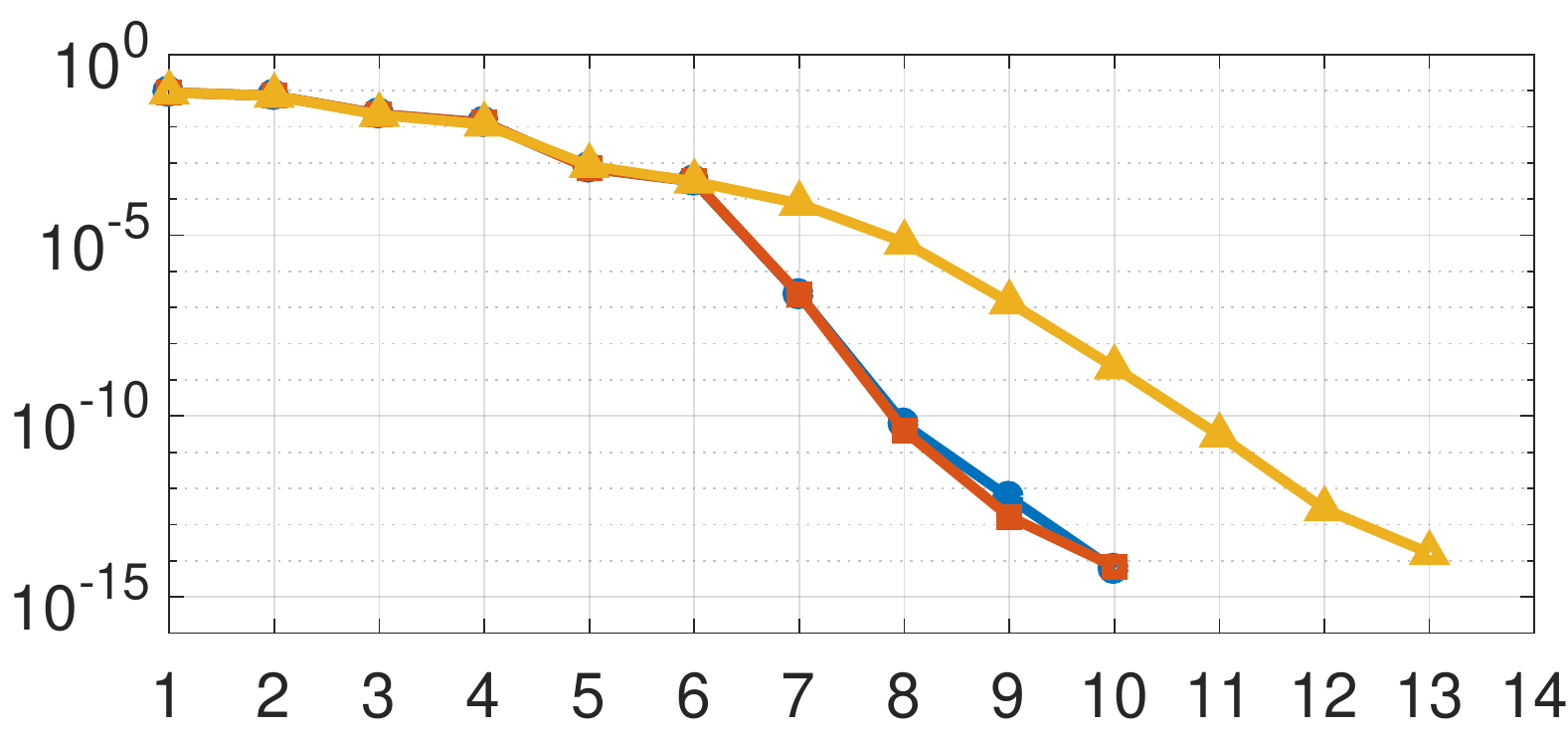}
		\hfill
		\includegraphics[height=2.8cm]{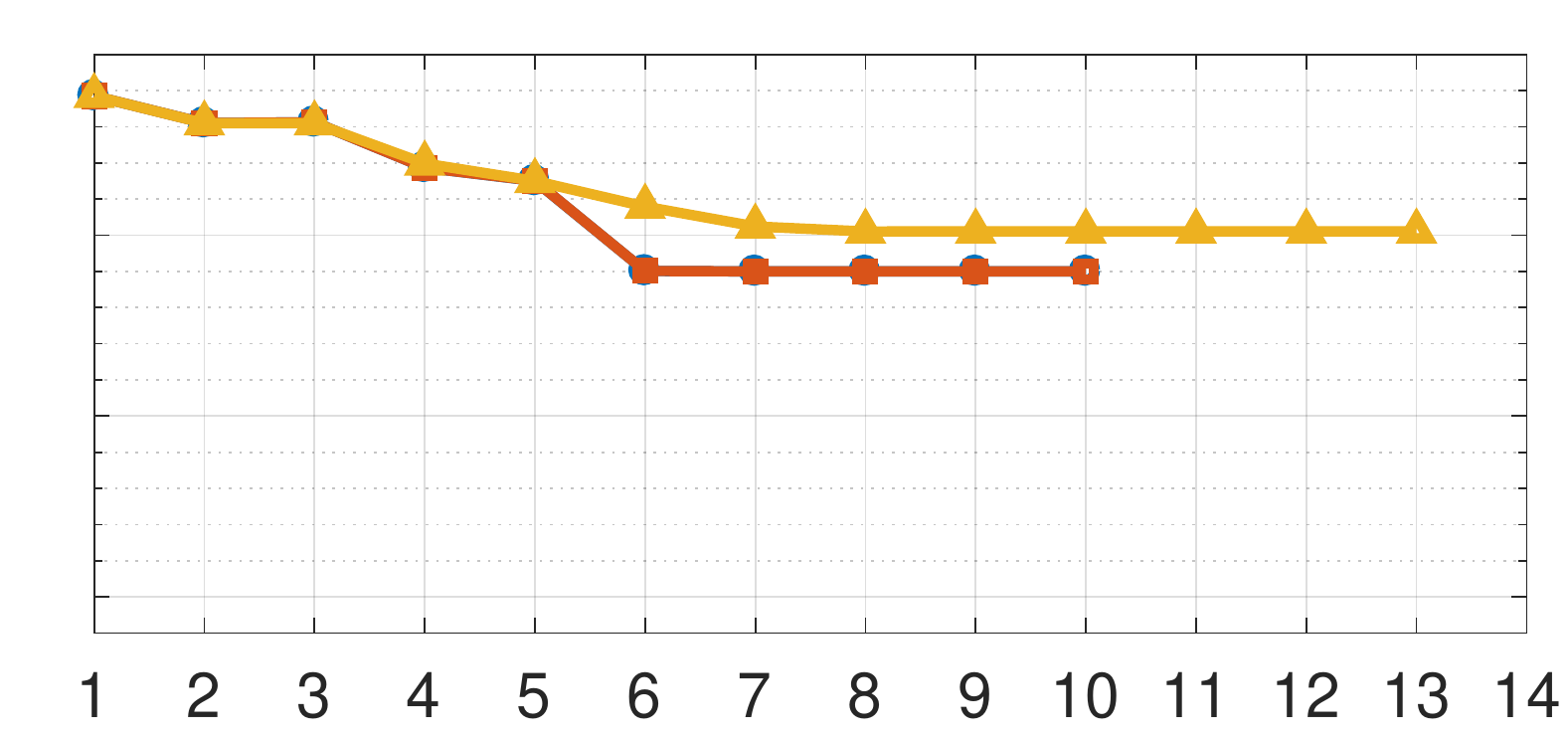}
	}
	
	\subfloat[Mesh 32x32]{
		\includegraphics[height=2.8cm]{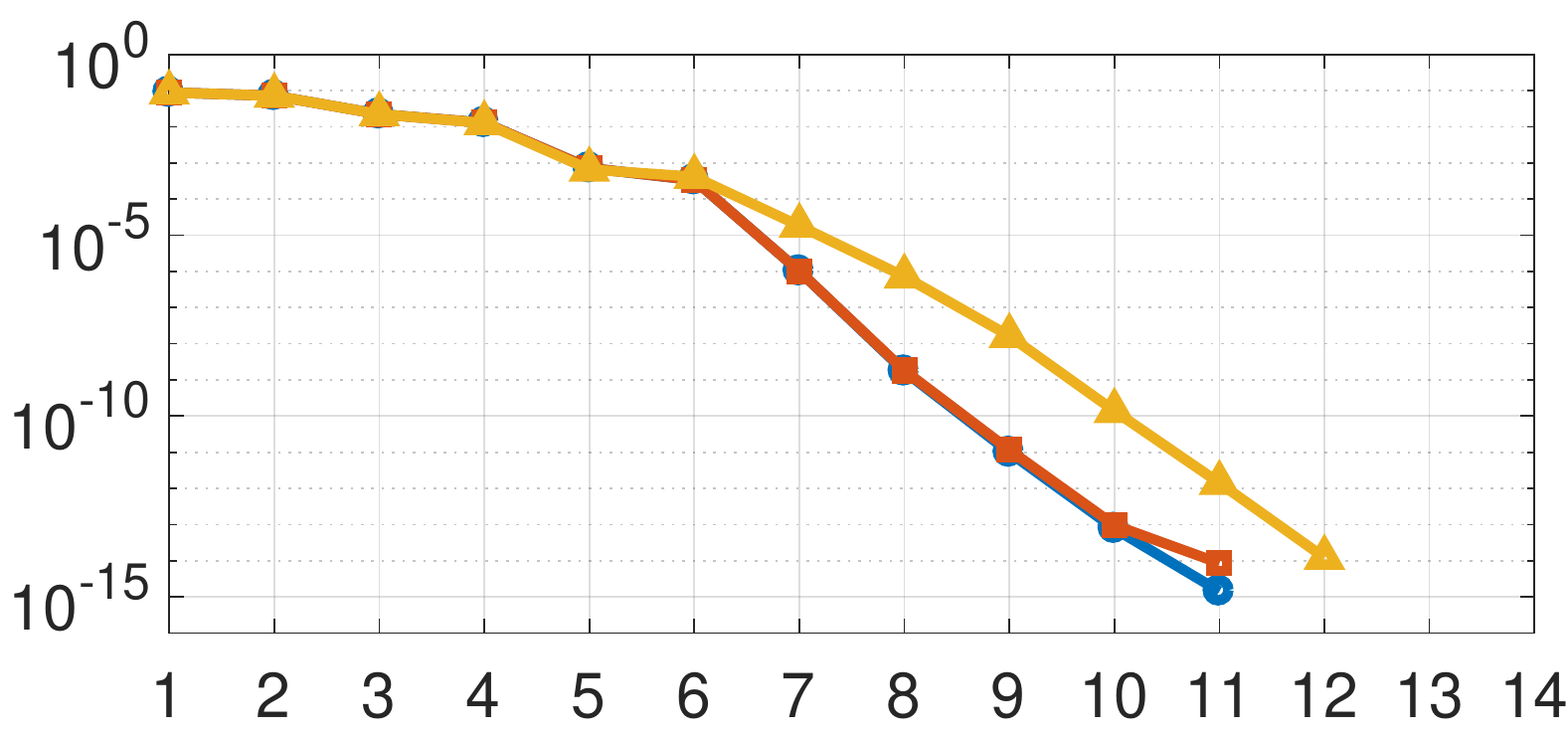}
		\hfill
		\includegraphics[height=2.8cm]{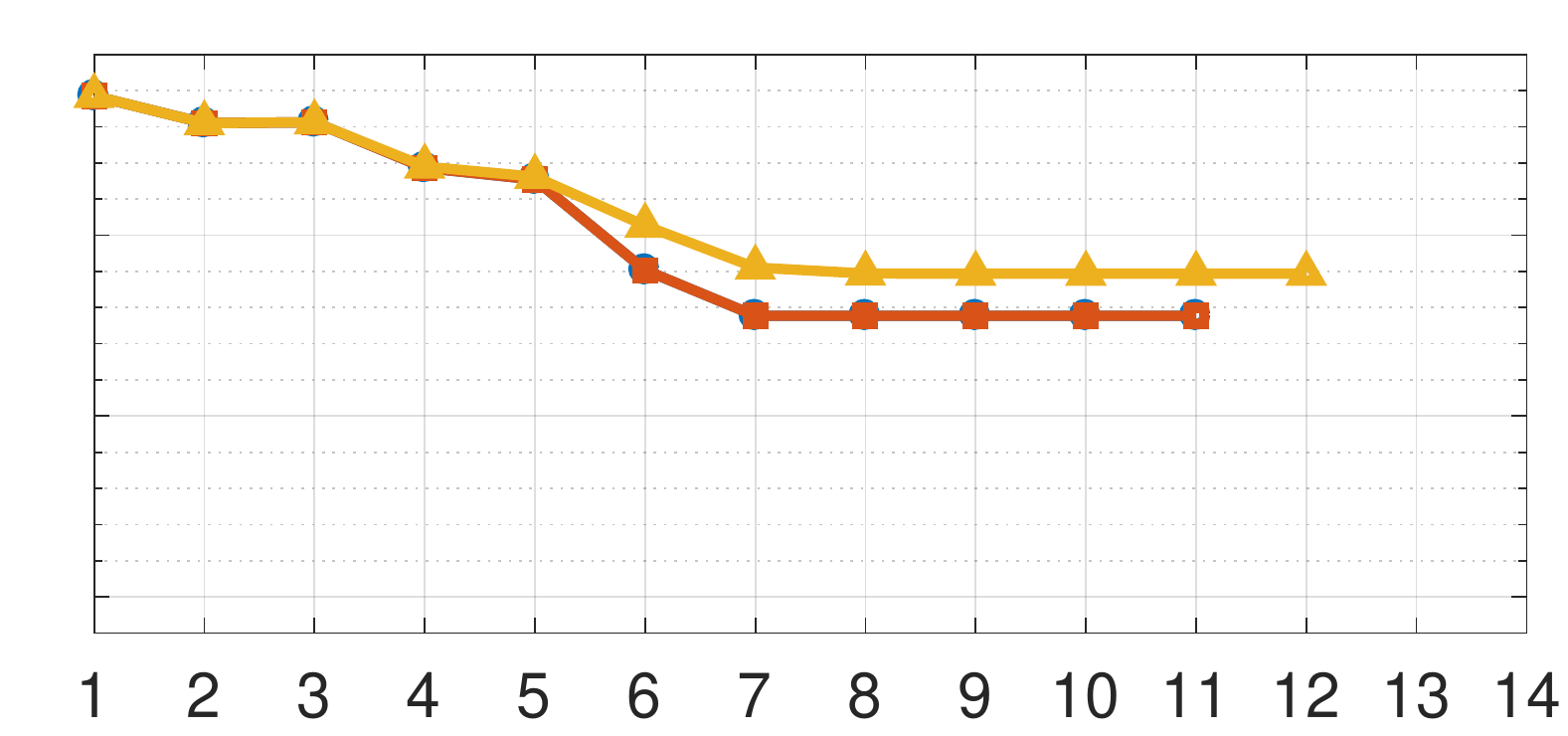}
	}
	
	\subfloat[Mesh 64x64]{
		\includegraphics[height=2.8cm]{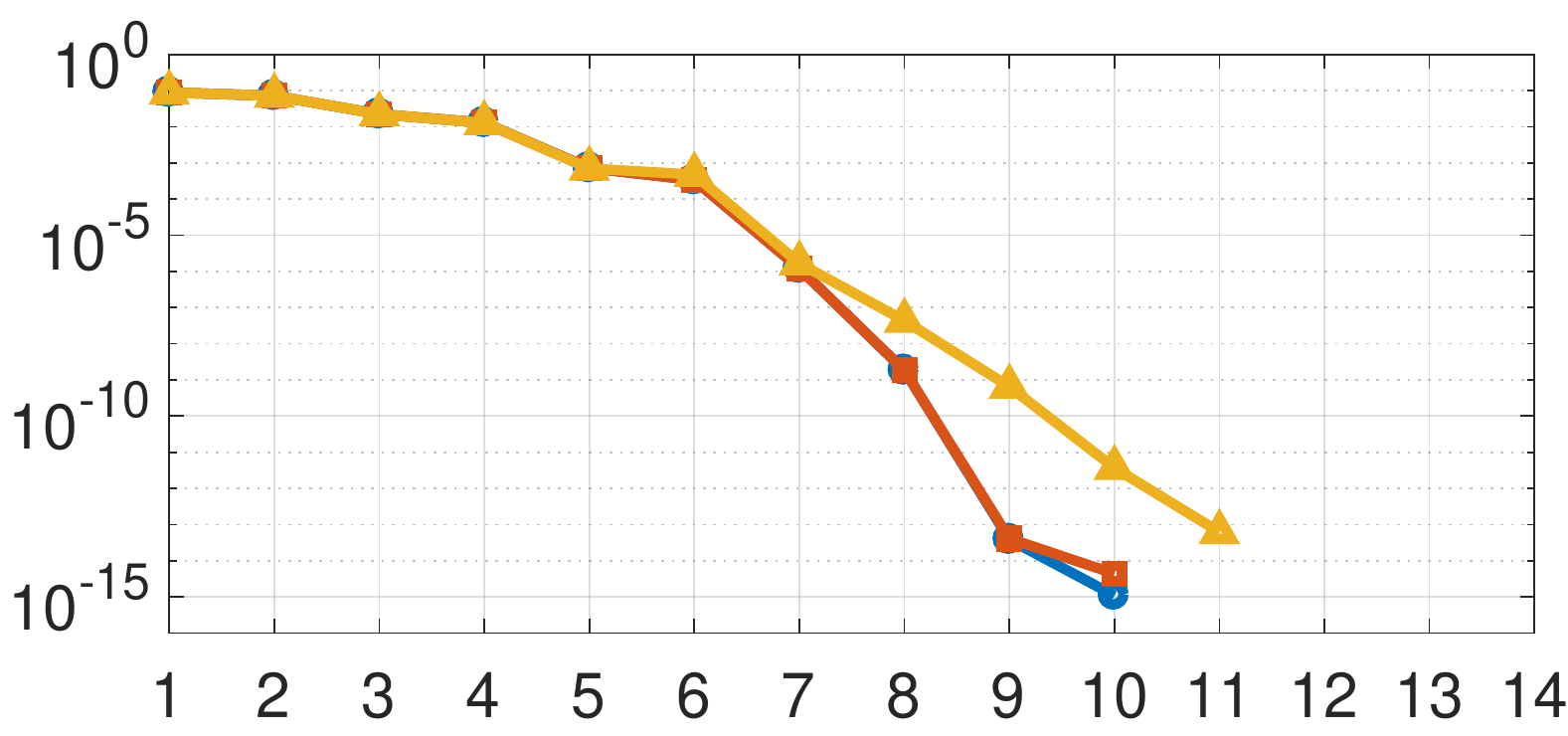}
		\hfill
		\includegraphics[height=2.8cm]{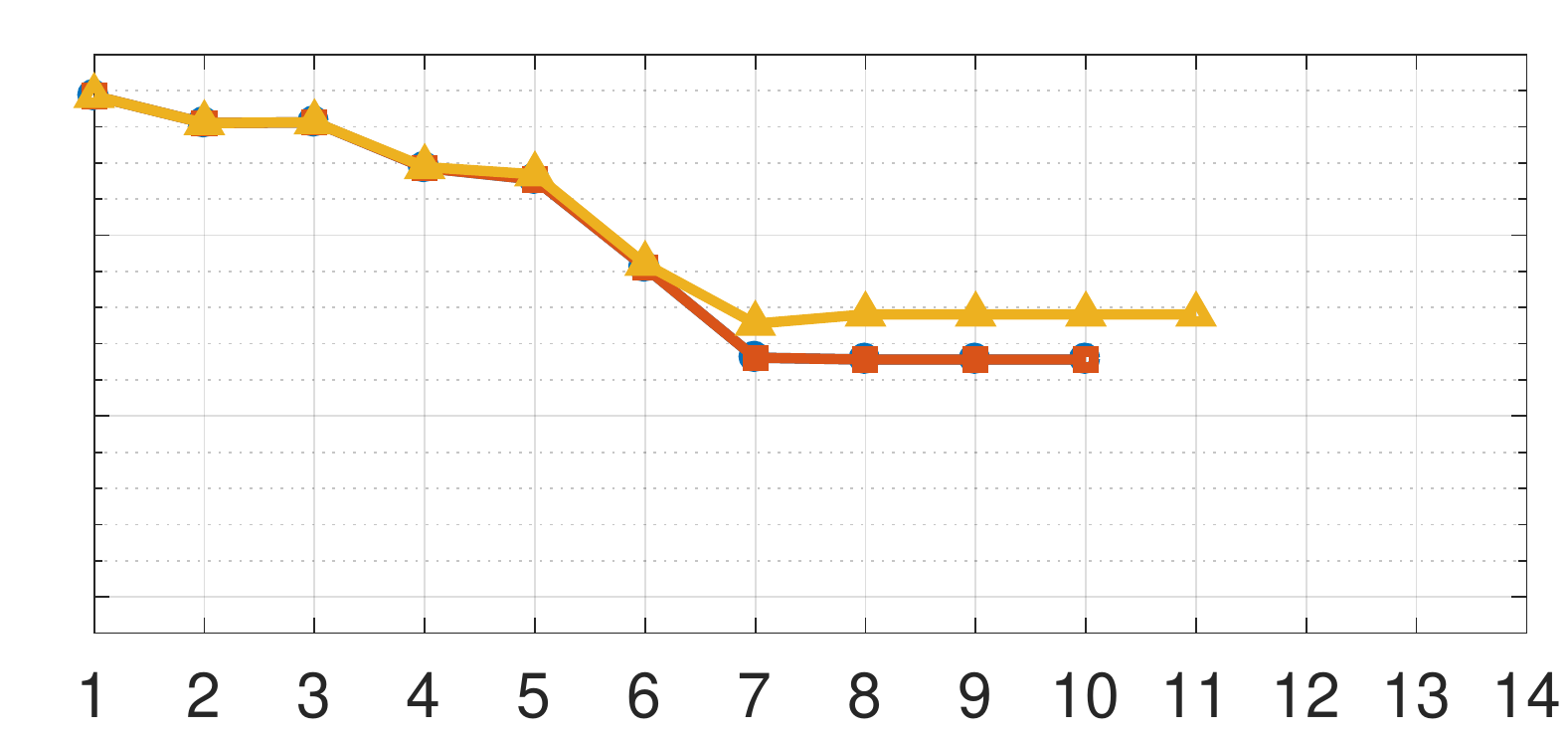}
	}

	\caption{A comparison of the three algorithms on Test 3 for different mesh sizes with cubic basis functions.}
	\label{fig:comparisonTest3}
\end{figure}

The introduction of the periodic conditions affects the behaviour of the three quasi-Newton schemes, but not dramatically.
As shown in Figure \ref{fig:comparisonTest3}, the algorithms require a couple extra iterations to reach the tolerance respect to the Dirichlet boundary condition case.
The convergence of the surface position error is also a bit rougher than in the previous cases.
However, the relative performances are not at all affected, and all three algorithms are still comparable.
As before Algorithms \ref{al:coupledSplitting} and \ref{al:decoupled} display essentially equal results.
In this test we kept the same choice for the initial guess for the free boundary: The flat curve $\Gamma_0 = \{ y = 1, \, 0 \leqslant x \leqslant 1 \}$.

Since the position of the exact free boundary does not lie in the trial functions space formed by the cubic B-splines basis, as in Test 2 a plateau is always reached, even though the level of the plateau is lowered with finer meshes.

\FloatBarrier


\section{Conclusions}
\label{sec:conclusions}

In this work we presented three different isogeometric-based algorithms for free boundary problems: Two follow a Galerkin approach and are an extension or modification of previously existing works, while one is a novel fully collocated scheme.
The dependence on the unknown geometry of the domain is handled through shape calculus, which results in a quasi-Newton method to be underlying the update strategy of the free boundary position.
While our interests in such algorithms is motivated by future applications, in the present paper we focused on giving a clear description of the implementation and numerical aspects.

We applied and compared the three algorithms to benchmark tests, with either Dirichlet or periodic boundary conditions on the lateral vertical sides of the domain.
The results show that, while having slight variations, the performances of all three algorithms are qualitatively comparable, and each of them converged to the correct solution of the problem.

The treatment of free boundary problems is computationally intense, especially in more complex problems.
For this reason the efficiency and speed of the algorithm is an important feature that needs to be taken into account.
In this respect, even if the collocated algorithm appeared to have slightly worse accuracy and sometimes required one or two extra iterations to reach the convergence tolerance, it proved to significantly outmatch the two Galerkin-based schemes on runtime, requiring in general less than half the time to complete the benchmarks.

Our future aim is now to apply the algorithms developed here to the resolution of the bifurcation branches of the Euler equations.
That problem presents several challenges due to the greater complexity of the equations and the intrinsic non-uniqueness of solutions at the bifurcation points, therefore both efficiency and precision are expected to play an important role.


\section{Acknowledgements}
\label{sec:acknowledgements}
MM and GS were partially supported by the European Research Council through the FP7 Ideas Consolidator Grant \emph{HIGEOM} n.616563.
FR was supported by grants no. 231668 and 250070 from the Norwegian Research Council.
This support is gratefully acknowledged.
MM and GS are members of the INdAM Research group GNCS.

\clearpage
\bibliographystyle{plain}
\bibliography{References}

\end{document}